\title{Conjugacy classes of reflections of maps}
\author{Gareth A. Jones\\
School of Mathematics\\
University of Southampton\\
Southampton SO17  1BJ, U.K.\\
{\tt G.A.Jones@maths.soton.ac.uk}
}

\documentclass[10pt]{article}
\usepackage{color}
\usepackage[latin1]{inputenc}
\usepackage{latexsym}
\usepackage{amsmath}
\usepackage{amssymb}
\usepackage{amsfonts}
\usepackage{tikz}
\newtheorem{thm}{Theorem}[section]
\newtheorem{lemma}[thm]{Lemma}
\newtheorem{cor}[thm]{Corollary}
\newtheorem{prop}[thm]{Proposition}
\date{}
\begin{document}

\maketitle

\begin{abstract}
This paper considers how many conjugacy classes of reflections a map can have, under various transitivity conditions. It is shown that for vertex- and for face-transitive maps there is no restriction on their number or size, whereas edge-transitive maps can have at most four classes of reflections. Examples are constructed, using topology, covering spaces and group theory, to show that various distributions of reflections can be achieved. Connections with real forms of algebraic curves are also discussed.
\end{abstract}

\medskip

\noindent{\bf MSC classification:} 05C10 (primary); 
14H37, 
14H57, 
20B25, 
30F10, 
30F50 
 (secondary).
 
 \medskip
 
 \noindent{\bf Key words:} Map, reflection, vertex-transitive, edge-transitive, conjugacy class, Riemann surface, real form.

\section{Introduction}

It is easily seen that if $\mathcal M$ is an orientably regular map then the number $cr({\mathcal M})$ of  conjugacy classes of reflections in its automorphism group ${\rm Aut}\,{\mathcal M}$ is bounded above by $3$, and that all values within this bound are attained (see Section~5.3). The aim of this note is to consider what can be said about conjugacy classes of reflections of maps if one relaxes the requirement of orientable regularity. 

In the case of vertex-transitive orientable maps, the following result shows that there are no group-theoretic  restrictions on the conjugacy classes of reflections, apart from the obvious fact that they must form a union of conjugacy classes of involutions in the non-trivial coset of a subgroup of index $2$ in the full automorphism group:

\begin{thm}\label{vtrans} Let $G$ be a finite group with a subgroup $G^+$ of index $2$, and let $K_1,\ldots, K_k$ be distinct conjugacy classes of involutions in $G\setminus G^+$ for some $k\ge 1$. Then there is a vertex-transitive map ${\mathcal M}$, on a compact orientable surface without boundary, such that ${\rm Aut}\,{\mathcal M}\cong G$ and the reflections of $\mathcal M$ correspond to the elements of the conjugacy classes $K_i$.
\end{thm}

It then follows that there are no arithmetic restrictions on the number of conjugacy classes of reflections, or on their sizes:
\begin{cor}\label{vtranscor}
If $c_1, \ldots, c_k$ are positive integers for some $k\ge 1$, there is a vertex-transitive map ${\mathcal M}$, on a compact orientable surface without boundary, such that the reflections of $\mathcal M$ form $k$ conjugacy classes of sizes $c_1, \ldots, c_k$. 
\end{cor}

Using map duality, one can replace the condition of vertex-transitivity in Theorem~\ref{vtrans} and Corollary~\ref{vtranscor} with face-transitivity. The situation is completely different for edge-transitive maps, as shown by the following result:

\begin{thm}\label{edgetrans}
If $\mathcal M$ is an edge-transitive map then ${\rm Aut}\,{\mathcal M}$ contains at most four conjugacy classes of reflections; if it has four then $\mathcal M$ is just-edge-transitive. 
\end{thm}

(An edge-transitive map is just-edge-transitive if it is neither vertex- nor face-transitive, or equivalently, if it has automorphism type~3 in the Graver-Watkins taxonomy of edge-transitive maps~\cite[Table~2]{GW}.) This result applies to all (connected) maps, including those which are non-compact, non-orientable, or with boundary. We will give examples (which can be chosen to be compact and without boundary) to show that any number $k\le 4$ of conjugacy classes of reflections can be realised by some edge-transitive map. In particular, in Theorem~5.1 we give an explicit construction of four infinite families of just-edge-transitive maps with $k=1, 2, 3$ or $4$ (the only possible values for such maps).

Algebraic geometry provides further motivation to study this topic. A compact Riemann surface $S$ is equivalent to a complex algebraic curve $C$, with the conjugacy classes of antiholomorphic involutions of $S$ corresponding to the real forms of $ C$, and the conjugacy classes of reflections (those with fixed points) corresponding to the real forms with real points (see~\cite{BCGG}, for example).
Bely\u\i's Theorem~\cite{Bel} shows that $C$ is defined over an algebraic number field if and only if the complex structure of $S$ is obtained in a standard way from a map $\mathcal M$ (or {\it dessin d'enfant}, in Grothendieck's terminology~\cite{Gro}). In this case the automorphisms of $\mathcal M$ are automorphisms (holomorphic or antiholomorphic) of $S$, and for many maps the converse is also true, so that the non-empty real forms of $C$ correspond to the conjugacy classes of reflections of $\mathcal M$.

This paper is organised as follows. Section~2 contains a mainly topological proof of Theorem~\ref{vtrans}, using Cayley graphs to construct the required vertex-transitive maps as regular coverings of single-vertex maps. Section~3 summarises a general approach to maps through permutation groups, used in Section~4 to study reflections.  These ideas are applied in Section~5 to edge-transitive maps, with a proof of Theorem~\ref{edgetrans} and the construction of several families of explicit examples. Section~6 discusses possible extensions to hypermaps. In Section~7 the results obtained here are compared with those obtained in the more general context of compact Riemann surfaces of a given genus by Bujalance, Gromadzki, Izquierdo, Natanzon and Singerman~\cite{BGI, BGS, GI, Nat}, with the work of Bujalance and Singerman~\cite{BS} on symmetry types of Riemann surfaces, and that of Meleko\u glu and Singerman~\cite{Mel, MS, MS14} on patterns of reflections of regular maps. Possible lines of future research are also discussed.

\section{Reflections of vertex-transitive maps}

In this paper we will consider maps on surfaces which may be orientable or non-orientable, with or without boundary, and compact or non-compact. They will always be assumed to be connected, and the main emphasis will be on compact orientable maps without boundary. By a {\em reflection\/} of a map we mean an automorphism (necessarily of order $2$) which fixes a point and acts as a reflection on a disc-shaped neighbourhood of that point. In the case of an orientable map, this implies that orientation is reversed, but the definition is also valid for non-orientable maps: for instance, the reflections of the icosahedron induce reflections of its antipodal quotient on the real projective plane.

We will say that a map $\mathcal M$ is {\em vertex}-, {\em edge}- or {\em face-transitive\/} if ${\rm Aut}\,{\mathcal M}$ acts transitively on its vertices, edges or faces; these conditions are all satisfied if $\mathcal M$ is {\em regular}, that is, if ${\rm Aut}\,{\mathcal M}$ acts transitively on its flags.

In order to prove Theorem~\ref{vtrans} and Corollary~\ref{vtranscor}, we first need the following:

\begin{lemma}\label{conjclasses}
Let $c_1, \ldots, c_k$ be positive integers for some $k\ge 1$. Then there is a finite group $G$ with a subgroup $G^+$ of index $2$ such that $G\setminus G^+$ contains a set of $k$ distinct conjugacy classes $K_i\;(i=1,\ldots, k)$, each consisting of $c_i$ involutions. 
\end{lemma}

\noindent{\sl Proof.} Let $G$ to be the direct product of $k$ dihedral groups $G_i$ of orders $2c_i$ or $4c_i$ as $c_i$ is odd or even; each $G_i$ contains the required class $K_i$ as a conjugacy class of reflections, and one can take $G^+$ to consist of the elements of $G$ with an even number of reflections among their coordinates.\hfill$\square$

\medskip

Note that in Lemma~\ref{conjclasses} we do not claim that the classes $K_i$ contain {\em all\/} the involutions in $G\setminus G^+$. For instance, if $k=1$ and $c_1=2$ then $K_1$ always generates a dihedral subgroup of $G$ contributing further involutions to $G\setminus G^+$. 

\medskip

\noindent{\sl Proof of Theorem~\ref{vtrans}.} Let $g_1,\ldots, g_k$ be representatives of the classes $K_i$. If they cannot be chosen to generate $G$, one can choose additional non-identity elements $g_{k+1},\ldots, g_l\in G^+$ so that $\langle g_1,\ldots, g_l\rangle=G$.

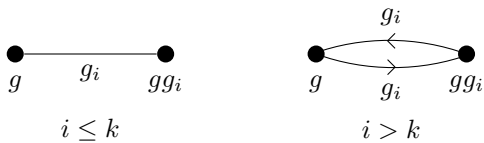
\begin{figure}[h!]
\label{involutions}
\begin{center}
 \begin{tikzpicture}[scale=0.5, inner sep=0.8mm]

\node (a) at (-8,0) [shape=circle, fill=black] {};
\node (b) at (-4,0) [shape=circle, fill=black] {};
\draw (a) to (b);

\node at (-8,-0.8) {$g$};
\node at (-4,-0.8) {$gg_i$};
\node at (-6,-0.5) {$g_i$};
\node at (-6,-2) {$i\le k$};


\node (c) at (0,0) [shape=circle, fill=black] {};
\node (d) at (4,0) [shape=circle, fill=black] {};
\draw [->] (4,0) arc (70:110:6);
\draw (0,0) arc (-110:-70:6);

\node at (0,-0.8) {$g$};
\node at (4,-0.8) {$gg_i$};

\draw (2.1,-0.35) to (1.9,-0.15);
\draw (2.1,-0.35) to (1.9,-0.55);
\draw (1.9,0.35) to (2.1,0.15);
\draw (1.9,0.35) to (2.1,0.55);
\node at (2,1) {$g_i$};
\node at (2,-1) {$g_i$};
\node at (2,-2) {$i>k$};

\end{tikzpicture}

\end{center}
\caption{Edges of $C$ corresponding to involutions $g_i$} 
\end{figure}

Now let $C$ be the Cayley graph for $G$ with respect to the generating set $X=\{g_1,\ldots, g_l\}$: this is a connected graph with vertex set $G$, and a directed edge labelled $g_i$ from $g$ to $gg_i$ for each $g\in G$ and each $g_i\in X$. The usual convention for Cayley graphs is that if a generator $g_i$ is an involution one replaces the directed edges from $g$ to $gg_i$ and from $gg_i$ to $g$ with a single undirected edge between these vertices; we will follow this rule if $i\le k$, but not if $i>k$ (see Figure~1), so that as an undirected graph $C$ has valency $m:=k+2(l-k)=2l-k$.

For each $h\in G$ the permutation $g\mapsto h^{-1}g$ of the vertices extends to an automorphism of $C$, preserving the directions and labelling of the edges. These automorphisms form a group $A\cong G$, acting regularly on the vertices of $C$. A simple calculation shows that the only elements of $G$ with fixed points on $C$ are those in the classes $K_i$: an element $h=gg_ig^{-1}\in K_i$, with $i\le k$, reverses the $|C_G(g_i)|/2=|G|/2c_i$ edges between pairs of vertices $g$ and $gg_i$, and fixes their midpoints. (If $g_i$ is an involution with $i>k$ then $h$ transposes the two directed edges between  $g$ and $gg_i$, so it has no fixed points.) The quotient graph $C'=A\backslash C$ therefore has a single vertex, incident with $k$ free edges labelled $g_1,\ldots, g_k$ and $l-k$ loops labelled $g_{k+1},\ldots, g_l$.

We now form a map $\mathcal M$ on a compact orientable surface $S$ which can be regarded as a tubular neighbourhood of $C$, with an action of $G$ on $\mathcal M$ imitating its action on $C$; however, we need to ensure that the involutions $g_i\;(i\le k)$ act on $S$ as reflections rather than half-turns. Let $T$ be an oriented surface formed from the $2$-sphere $S^2$ by removing $m$ open discs, leaving mutually disjoint boundary components $B_i\;(i=1,\ldots, m)$, each homeomorphic to $S^1$. Intuitively, if $i^*$ is defined to be $i$ or $i+l-k$ as $i\le k$ or $i>k$, one can think of $T$ as surrounding a typical vertex of $C$, with undirected edges labelled $g_i$ leaving though $B_i$ for $i=1,\ldots, k$, and pairs of directed edges labelled $g_i$ leaving through $B_i$ and entering through $B_{i*}$ for $i=k+1,\ldots, l$. We construct a map $\mathcal T$ on $T$ by choosing a vertex $v$ in the interior of $T$, drawing disjoint free edges $e_i$ joining $v$ to a point $p_i$ (not a vertex) in each $B_i$, and then drawing mutually disjoint loops $l_i$ in $T$ from $v$ to $v$, enclosing $e_i$ and $B_i$ but not $e_j$ or $B_j$ for any $j\ne i$ (see Figure~2). This map has $m+1$ faces: there are $m$ $4$-gons $f_i$, each bounded by $e_i$ (twice), $l_i$ and $B_i$ for some $i=1,\ldots, m$, and there is one $m$-gon $f$, bounded by $l_1,\ldots, l_m$. (The closure of each face $f_i$ is homeomorphic to a rectangle with opposite sides identified, that is, a cylinder.)

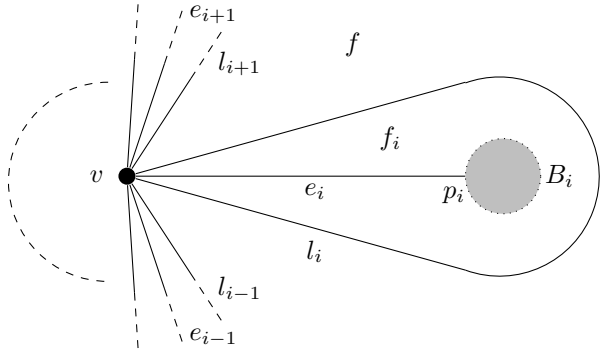
\begin{figure}[h!]
\label{K5}
\begin{center}
 \begin{tikzpicture}[scale=0.5, inner sep=0.8mm]
 
\node (c) at (0,0) [shape=circle, fill=black] {};
\node at (9,0) {};
\draw (c) to (9,0); \draw  [dotted] [fill=lightgray] (11,0) arc (0:360:1);
 
\draw (c) to (1,3);
\draw (c) to (1.7,2.6);
\draw (c) to (0.2,3.1);

\draw [dashed] (1,3) to (1.5,4.5);
\draw [dashed] (1.7,2.6) to (2.55,3.9);
\draw [dashed] (0.2,3.1) to (0.3,4.65);

\node at (2.3,4.3) {$e_{i+1}$};
\node at (3,3) {$l_{i+1}$};

\draw (c) to (1,-3);
\draw (c) to (1.7,-2.6);
\draw (c) to (0.2,-3.1);

\draw [dashed] (1,-3) to (1.5,-4.5);
\draw [dashed] (1.7,-2.6) to (2.55,-3.9);
\draw [dashed] (0.2,-3.1) to (0.3,-4.65);

\node at (2.3,-4.3) {$e_{i-1}$};
\node at (3,-3) {$l_{i-1}$};

\draw (c) to (9,2.5);
\draw (c) to (9,-2.5);
\draw (9,-2.5) arc (-110:110:2.65);
 
  \draw [dashed] (-0.5,2.5) arc (90:270:2.65);

\node at (-0.8,0) {$v$};
\node at (5,-0.4) {$e_i$};
\node at (5,-2) {$l_i$};
\node at (11.5,0) {$B_i$};
\node at (8.7,-0.5) {$p_i$};
\node at (7,1) {$f_i$};
\node at (6,3.5) {$f$};

\end{tikzpicture}

\end{center}
\caption{The map $\mathcal T$ on $T$} 
\end{figure}
We construct $\mathcal M$ by joining $|G|$ copies of the map $\mathcal T$ across their boundary components. More precisely, we form a quotient $S$ of $T\times G$ by identifying each $(p,g)\in B_i\times\{g\}$ with $(p,gg_i)\in B_i\times\{gg_i\}$ for $i=1,\ldots,k$, and identifying $B_i\times\{g\}$ with $B_{i^*}\times\{gg_i\}$ for $i=k+1,\ldots, l$ by means of a homeomorphism $B_i\to B_{i^*}$ which sends $p_i$ to $p_{i^*}$ and reverses the orientations of these two boundary components induced from that of $T$. Thus the orientation of $T$ is reversed or preserved when passing between $T\times\{g\}$ and $T\times\{gg_i\}$ for $i\le k$  or for $i>k$ respectively; since the generators $g_i$ lie in $G\setminus G^+$ or $G^+$ as $i\le k$ or $i>k$, this means that the resulting surface $S$ is orientable.

The map $\mathcal T$ induces isomorphic maps on each $T\times\{g\}$ and hence, via these identifications, it induces a map $\mathcal M$ on $S$. This has $|G|$ vertices $(v,g)$ where $g\in G$. Each vertex $(v,g)$ is incident with $m$ loops $(l_i,g)$ for  $i=1,\ldots, m$, with $m$ edges of the form $(e_i,g)\cup(e_{i^*}, gg_i)$ for $i=1,\ldots,m$, and with $l-k$ edges $(e_i,gg_i^{-1})\cup(e_{i^*},g)$ for $i=k+1,\ldots, l$. There are $|G|$ $m$-gonal faces $(f,g)$, where $g\in G$, and $|G|m/2$ $4$-gonal faces $(f_i,g)\cup(f_{i^*},gg_i)$ where $g\in G$ and $i=1,\ldots, m$; the latter are bounded by $(l_i,g)$, $(l_{i^*},gg_i)$ and $(e_i,g)\cup(e_{i^*}, gg_i)$ (twice).

Now let $G$ act on $T\times G$, with each $h\in G$ sending $(p,g)$ to $(p,h^{-1}g)$ for each point $p\in T$ and $g\in G$, so that $G$ permutes the $|G|$ copies $T\times\{g\}$ of $T$ regularly. This action is compatible with the identifications, and therefore induces a faithful action of $G$ as a vertex-transitive group of automorphisms of $\mathcal M$. By the construction, the only points in $S$ fixed by a non-identity element of $G$ are those in the identified boundary components $B_i\times \{g\} = B_{i^*}\times\{gg_i\} = B_i\times\{gg_i\}$ for $g\in G$ and $i=1,\ldots, k$: these are fixed by the involution $gg_ig^{-1}\in K_i$, acting as a reflection of $\mathcal M$. We have $G\le{\rm Aut}\,{\mathcal M}$, and we can ensure that these groups are equal by modifying $\mathcal T$ (and hence $\mathcal M$) with the addition of additional loops attached to $v$ within the face $f$. This guarantees that $\mathcal M$ has no reflections other than those in the classes $K_i$.  \hfill$\square$

\medskip

\noindent{\sl Proof of Corollary~\ref{vtranscor}.} This follows from Lemma~\ref{conjclasses} and Theorem~\ref{vtrans}. \hfill$\square$

\medskip

In the construction used in the proof of Theorem~\ref{vtrans}, $\mathcal M$ has $|G|$ vertices, $3|G|m/2$ edges, and $|G|(1+m/2)$ faces, so it has Euler characteristic
\[|G|\Bigl(1-\frac{3m}{2}+1+\frac{m}{2}\Bigr)=|G|(2-m).\]

\medskip

\noindent{\bf Example 2.1} Let $G=S_n$ and $k=1$, with $K_1$ the class of transpositions (one of $\lfloor(n+2)/4\rfloor$ conjugacy classes of odd involutions). We can take $g_1=(1,2)$ as a representative of $K_1$, adjoining $g_2=(1, 2, \ldots, n)$ or $(1, 2)(1, 2, \ldots, n)=(1, 3, 4, \ldots, n)$ as $n$ is odd or even to form a generating set for $G$. Thus $l=2$ and $m=3$, so the resulting map $\mathcal M$ has characteristic $-n!$ and genus $1+\frac{1}{2}n!$.

If $n=4$, for instance, then $C$ is the $1$-skeleton of a truncated cube. Each of the six transpositions $h\in K_1$ acts on $C$ as a half-turn of the cube, fixing the midpoints of two antipodal edges. This illustrates why, in the construction of $\mathcal M$, we needed to be careful about the identifications of the boundary components. If we had simply formed $S$ hy using a tubular neighbourhood of $C$ in ${\mathbb R}^3$, then each $h$ would act as a half-turn on $S$, fixing four points. Instead, by our identifications of the boundary components $B_i\times\{g\}$ for $i\le k$, equivalent to cutting the corresponding tubes and rejoining them with reverse orientation, we ensured than $h$ acts as a reflection of $S$, fixing the whole of $B_i\times\{g\}$ rather than just two of its points.

\section{Algebraic map theory}

In order to make further progress in the area, and in particular to consider possible generalisations of Theorem~\ref{vtrans} and Corollary~\ref{vtranscor}, we need an alternative approach to maps, based on group theory. In this section we will briefly outline the algebraic theory of maps developed in more detail elsewhere~\cite{BS, JS}. 

Maps on surfaces $S$ (possibly non-orientable and possibly with boundary) correspond to permutation representations of the group
\[\Gamma=\langle R_0, R_1, R_2\mid R_i^2=(R_0R_2)^2=1\rangle,\]
the free product of a Klein four-group $E=\langle R_0, R_2\rangle$ and a cyclic group $\langle R_1\rangle$ of order~$2$. Given a map $\mathcal M$, let $\Phi$ denote the set consisting of its flags $\phi=(v,e,f)$, where $v, e$ and $f$ denote a mutually incident vertex, edge and face of $\mathcal M$. For each such $\phi$ and each $i=0, 1, 2$, there is at most one flag $\phi'\ne \phi$ sharing the same $j$-dimensional components as $\phi$ for each $j\ne i$ (there may be none if $\phi$ is a boundary flag). Define $r_i$ to be the permutation of $\Phi$ which transposes each $\phi$ with $\phi'$ if the latter exists, and fixes $\phi$ if it does not (see Figure~3 for the former case). Then
\[r_i^2=(r_0r_2)^2=1,\]
so there is a permutation representation
\[\theta:\Gamma\to G:=\langle r_0, r_1, r_2\rangle\le{\rm Sym}\,\Phi\]
of $\Gamma$ on $\Phi$, given by $R_i\mapsto r_i$.

\begin{figure}[h!]
\label{K5}
\begin{center}
 \begin{tikzpicture}[scale=0.5, inner sep=0.8mm]

\node (c) at (0,0) [shape=circle, fill=black] {};
\node (d) at (8,0) [shape=circle, fill=black] {};
\draw (c) to (d);
\draw (c) to (1,-3);
\draw (c) to (1,3);
\draw (d) to (7,-3);
\draw (d) to (7,3);
\draw (c) to (-2.5,2.5);
\draw (c) to (-2.5,-2.5);
\draw (d) to (10.5,2.5);
\draw (d) to (10.5,-2.5);

\draw (c) to (1,0.5);
\draw (c) to (1,-0.5);
\draw (1,0.5) to (1,-0.5);
\draw (d) to (7,0.5);
\draw (d) to (7,-0.5);
\draw (7,0.5) to (7,-0.5);
\draw (c) to (0.8,0.8);
\draw (0.3,1) to (0.8,0.8);

\node at (-0.8,0) {$v$};
\node at (4,-0.4) {$e$};
\node at (4,2) {$f$};

\node at (1.5,0.5) {$\phi$};
\node at (6.2,0.5) {$\phi r_0$};
\node at (1.4,1.3) {$\phi r_1$};
\node at (1.9,-0.5) {$\phi r_2$};
\node at (5.9,-0.5) {$\phi r_0r_2$};

\end{tikzpicture}

\end{center}
\caption{Generators $r_i$ acting on a flag $\phi=(v,e,f)$} 
\end{figure}
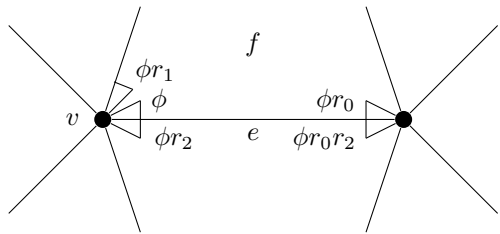

Conversely, given a permutation representation of $\Gamma$ on any set $\Phi$, one can construct a map $\mathcal M$ by identifying the vertices, edges and faces of $\mathcal M$ with the orbits of the dihedral subgroups $\langle R_1, R_2\rangle\cong D_{\infty}$,  $\langle R_0, R_2\rangle\cong D_2\cong V_4$ and  $\langle R_0, R_1\rangle\cong D_{\infty}$ on $\Phi$, with incidence given by non-empty intersection.

The map $\mathcal M$ is connected if and only if $\Gamma$ acts transitively on $\Phi$ (as we will always assume), in which case the stabilisers in $\Gamma$ of flags $\phi\in\Phi$ form a conjugacy class of subgroups $M\le\Gamma$, known as {\em map subgroups}. 
One can regard $\Gamma$ as the automorphism group of the {\em universal map\/} ${\mathcal M}_{\infty}$ (see Figure~4), so that $\mathcal M$ is isomorphic to the quotient ${\mathcal M}_{\infty}/M$ of ${\mathcal M}_{\infty}$ by a map subgroup $M\le\Gamma$. 

\begin{figure}[h!]
\label{univbipmap}
\begin{center}
 \begin{tikzpicture}[scale=0.5, inner sep=0.5mm]

\node (A) at (-10,0) [shape=circle, fill=black, label=below: $\frac{-1}{1}$] {};
\node (B) at (0,0) [shape=circle, fill=black, label=below: $\frac{0}{1}$] {};
\node (C) at (10,0) [shape=circle, fill=black, label=below: $\frac{1}{1}$] {};
\node (D) at (-6.7,0) [shape=circle, fill=black, label=below: $\frac{-2}{3}$] {};
\node (E) at (-3.3,0) [shape=circle, fill=black, label=below: $\frac{-1}{3}$] {};
\node (F) at (3.3,0) [shape=circle, fill=black, label=below: $\frac{1}{3}$] {};
\node (G) at (6.7,0) [shape=circle, fill=black, label=below: $\frac{2}{3}$] {};
\node (H) at (-8,0) [shape=circle, fill=black, label=below: $\frac{-4}{5}$] {};
\node (I) at (-6,0) [shape=circle, fill=black, label=below: $ $] {};
\node (J) at (-4,0) [shape=circle, fill=black, label=below: $ $] {};
\node (K) at (-2,0) [shape=circle, fill=black, label=below: $\frac{-1}{5}$] {};
\node (L) at (2,0) [shape=circle, fill=black, label=below: $\frac{1}{5}$] {};
\node (M) at (4,0) [shape=circle, fill=black, label=below: $\frac{2}{5}$] {};
\node (N) at (6,0) [shape=circle, fill=black, label=below: $\frac{3}{5}$] {};
\node (O) at (8,0) [shape=circle, fill=black, label=below: $\frac{4}{5}$] {};

\draw (B) arc (0:180:5);
\draw (C) arc (0:180:5);
\draw (D) arc (0:180:1.65);
\draw (B) arc (0:180:1.65);
\draw (F) arc (0:180:1.65);
\draw (C) arc (0:180:1.65);
\draw (H) arc (0:180:1);
\draw (B) arc (0:180:1);
\draw (L) arc (0:180:1);
\draw (C) arc (0:180:1);
\draw (I) arc (0:180:0.35);
\draw (E) arc (0:180:0.35);
\draw (M) arc (0:180:0.35);
\draw (G) arc (0:180:0.35);

\draw[dotted] (-11,0) to (11,0);

\node at (12,0) {$\mathbb Q$};
\node at (12,4) {$\mathbb H$};

\draw [dashed] (B) to (0,7);
\draw [dashed] (5,0) to (5,7);
\draw [<->] (-0.8,6) to (0.8,6);
\draw [<->] (4.2,6) to (5.8,6);
\draw [<->] (1.6,4.6) to (2.55,3.45);

\node at (6.5,6) {$R_0$};
\node at (-1.5,6) {$R_1$};
\node at (3,3) {$R_2$};

\end{tikzpicture}

\end{center}
 \caption{The universal map ${\mathcal M}_{\infty}$ with reflections $R_i$} 
\end{figure}
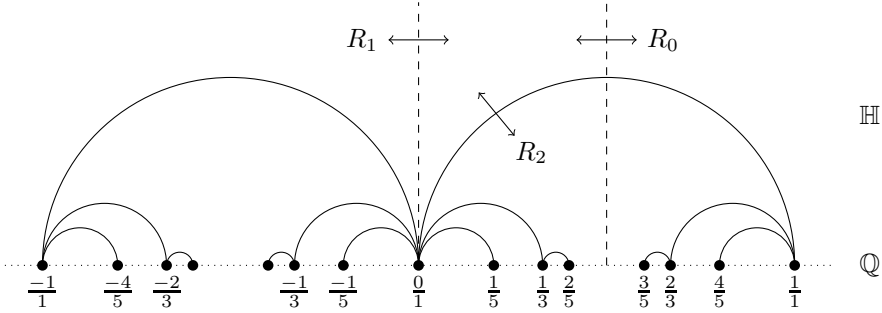

The underlying surface of ${\mathcal M}_{\infty}$ is the upper half plane ${\mathbb H}=\{z\in{\mathbb C}\mid {\rm Im}\,z>0\}$, with vertices at the reduced rationals $a/b$ with $b$ odd, and an edge (a hyperbolic geodesic in $\mathbb H$) between vertices $a/b$ and $c/d$ if and only if $ad-bc=\pm 1$. Part of this map, with $|{\rm Re}\,z|\le 1$ and $0<b\le 5$, is shown in Figure~4; the pattern repeats with period $1$. Then $\Gamma$ is faithfully represented as ${\rm Aut}\,{\mathcal M}_{\infty}$,  the principal congruence subgroup of level $2$ in $PGL_2({\mathbb Z})$, with its generators $R_i$ acting as the reflections
\[R_0:z\mapsto 1-\overline z,\quad R_1: z\mapsto-\overline z,\quad R_2: z\mapsto 
\frac{\overline z}{2\overline z-1}.\]

The map $\mathcal M$ is compact if and only if $\Phi$ is finite, and it has non-empty boundary if and only if some $r_i$ has fixed points in $\Phi$. More generally, $\mathcal M$ is orientable and without boundary if and only if $M$ is contained in the even subgroup $\Gamma^+$ of $\Gamma$, consisting of the words of even length in the generators $R_i$.

The group $G$ is known as the {\em monodromy group\/} ${\rm Mon}\,{\mathcal M}$ of $\mathcal M$. The {\em automorphism group\/} $A={\rm Aut}\,{\mathcal M}$ of $\mathcal M$ is the centraliser of $G$ in ${\rm Sym}\,\Phi$, isomorphic to $N/M$ where $N$ is the normaliser $N_{\Gamma}(M)$ of $M$ in $\Gamma$. The map $\mathcal M$ is {\em regular\/} if $A$ acts transitively on $\Phi$; this is equivalent to $G$ being a regular permutation group, that is, to $M$ being normal in $\Gamma$, in which case
\[A\cong G\cong \Gamma/M.\]

The reflections $R$ of ${\mathcal M}_{\infty}$ are the conjugates of the generating reflections $R_i$. The reflections $r$ of $\mathcal M$ are its automorphisms induced by reflections $R$ of ${\mathcal M}_{\infty}$ contained in $N$ (see~\cite[Lemma 1.5.9]{BCGG} for the analogous result for Riemann surfaces). If two reflections $R, R'\in N$ are conjugate in $N$ then their images $r, r'$ are conjugate in ${\rm Aut}\,{\mathcal M}$. It follows that the number $cr(N)$ of conjugacy classes of reflections in $N$ is an upper bound for the number $cr({\mathcal M})$ of conjugacy classes of reflections of $\mathcal M$, and that $cr({\mathcal M})\ge 1$ if $cr(N)\ge 1$.

The barycentric subdivision $B({\mathcal M})$ of a map $\mathcal M$ is a triangulation of $S$, with its faces corresponding to the flags of $\mathcal M$. The dual map $B({\mathcal M})^*$, which has its vertices corresponding to the flags of $\mathcal M$, can be regarded as an embedding of a permutation graph for $\Gamma$ on $\Phi$, or equivalently a Schreier coset graph $\Sigma$ for $M$ in $\Gamma$, with respect to the generators $R_i$ of $\Gamma$: there is an undirected edge, labelled $i$, between pairs of vertices (or cosets of $M$) if and only if they form a $2$-cycle of $R_i$ (we can and will omit loops at vertices corresponding to fixed points). 

Any spanning tree $T$ for this graph $\Sigma$ yields a Schreier transversal for $M$ in $\Gamma$: an arbitrary vertex is chosen as a base point, representing $M$, and then the labels of successive edges in the unique shortest path in $T$ to any other vertex give a word in the generators $R_i$ which serves as a representative for the corresponding coset. The Reidemeister-Schreier algorithm~\cite[\S II.4]{LS}, \cite[\S 2.3]{MKS}, applied to this transversal, then gives a presentation for $M$. As a particular case of the Kurosh Subgroup Theorem~\cite[\S IV.1]{LS}, \cite[\S 4.3]{MKS}, we find that $M$ is a free product of subgroups conjugate in $\Gamma$ to $\langle R_1\rangle$ or to a non-identity subgroup $F$ of $E=\langle R_0, R_2\rangle$, and a free group. The finite factors in this decomposition are representatives of the conjugacy classes of maximal finite subgroups of $M$, corresponding to flags on the boundary of $\mathcal M$, fixed by $R_1$ or by $F$; the free factor, isomorphic to the quotient of $M$ by the normal subgroup generated by its torsion elements, can be identified with the fundamental group of the surface $S_0$ formed by puncturing $S$ at the vertices and face-centres of $\mathcal M$. Indeed, maps can  be regarded as orbifolds, so that $M$ itself is the orbifold fundamental group of $\mathcal M$, regarded as a covering of the trivial map or orbifold ${\mathcal M}_{\infty}/\Gamma$. 

\medskip

\noindent{\bf Example 3.1} Let $\mathcal M$ be the map on the closed unit disc $D=\{z\in{\mathbb C}\mid |z|\le 1\}$ with two vertices at $\pm 1$ joined by an edge along $\mathbb R$ (see Figure~5). This is a regular map with four flags, fixed by $R_1$ and permuted regularly by $E$, so $M$ is the normal closure of $R_1$ in $\Gamma$. By taking a Schreier transversal $1$, $R_2$, $R_0R_2$, $R_0$ for $M$ in $\Gamma$ we find that $M$ has generators $S_i=R_1$, $R_1^{R_2}$, $R_1^{R_0R_2}$ and $R_1^{R_0}$ for $i=1,\ldots, 4$, with defining relations $S_i^2=1$, so $M$ is the free product of four copies of $C_2$. The punctured surface $D_0$, a closed disc minus four boundary points, is simply connected, so the free part of the decomposition of $M$ is trivial.

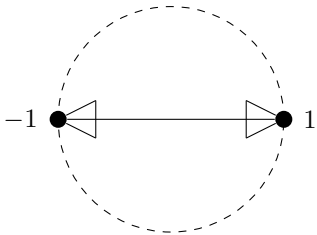
\begin{figure}[h!]
\label{4flagmap}
\begin{center}
 \begin{tikzpicture}[scale=0.5, inner sep=0.8mm]

\node (c) at (0,0) [shape=circle, fill=black] {};
\node (d) at (6,0) [shape=circle, fill=black] {};
\draw (c) to (d);

\draw (c) to (1,0.5);
\draw (c) to (1,-0.5);
\draw (1,0.5) to (1,-0.5);
\draw (d) to (5,0.5);
\draw (d) to (5,-0.5);
\draw (5,0.5) to (5,-0.5);
\draw [dashed] (6,0) arc (0:360:3);

\node at (-1,0) {$-1$};
\node at (6.7,0) {$1$};

\end{tikzpicture}

\end{center}
\caption{A map on the closed disc} 
\end{figure}

In studying reflections $r$ of a map $\mathcal M$, it is useful to consider those flags $\phi=(v,e,f)$ of $\mathcal M$ which have a component incident with the subset ${\rm Fix}(r)$ of $\mathcal M$ fixed by $r$, that is, which have at least one component (in fact, always two) invariant under $r$. This means that $r$ sends $\phi$ to one of its images $\phi r_i$ under the standard generators $R_i\;(i=0, 1, 2)$ of $\Gamma$, so we will say that $r$ acts with type $i$ on $\phi$; this corresponds to some conjugate of $R_i$ in $N:=N_{\Gamma}(M)$ inducing the automorphism $r\in{\rm Aut}\,{\mathcal M}\cong N/M$ (see Figure~6).

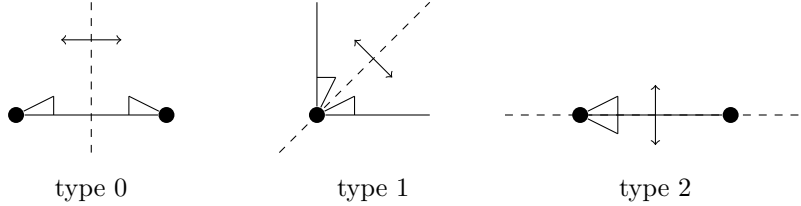
\begin{figure}[h!]
\label{K5}
\begin{center}
 \begin{tikzpicture}[scale=0.5, inner sep=0.75mm]

\node (a) at (-8,0) [shape=circle, fill=black] {};
\node (b) at (-4,0) [shape=circle, fill=black] {};
\draw (a) to (b);
\draw (-7,0.5) to (a);
\draw (-7,0.5) to (-7,0);
\draw (-5,0.5) to (b);
\draw (-5,0.5) to (-5,0);
\draw [dashed] (-6,3) to (-6,-1);
\draw [<->] (-6.8,2) to (-5.2,2);
\node at (-6,-2) {type $0$};


\node (c) at (0,0) [shape=circle, fill=black] {};
\draw (c) to (3,0);
\draw (c) to (0,3);
\draw (1,0.5) to (c);
\draw (1,0.5) to (1,0);
\draw (0.5,1) to (c);
\draw (0.5,1) to (0,1);
\draw [dashed] (-1,-1) to (3,3);
\draw [<->] (2,1) to (1,2);
\node at (1.5,-2) {type $1$};


\node (d) at (7,0) [shape=circle, fill=black] {};
\node (e) at (11,0) [shape=circle, fill=black] {};
\draw (d) to (e);
\draw (8,0.5) to (d);
\draw (8,-0.5) to (d);
\draw (8,0.5) to (8,-0.5);
\draw [dashed] (5,0) to (13,0);
\draw [<->] (9,0.8) to (9,-0.8);
\node at (9,-2) {type $2$};

\end{tikzpicture}

\end{center}
\caption{Actions of types $0, 1, 2$ on flags} 
\end{figure}

Of course, $r$ may have actions of different types on different flags: for instance, if $r$ has an action of type $2$ on $\phi$, fixing $v$ and $e$, and $v$ has odd valency $n$, then $r$ has an action of type $1$ on the flag $\phi(r_1r_2)^{(n-1)/2}$ incident with $v$. There is a similar pairing of actions of types $0$ and $1$ arising from faces of odd valencies. Even if all valencies are even, so that the type is constant along each connected component of ${\rm Fix}(r)$, it is possible for different components to yield actions of different types: for instance, if $b$ is odd then the torus map $\{4,4\}_{b,0}$ has reflections $r$ (in axes parallel to edges) with two components of fixed points, each yielding actions entirely of type $0$ or of type $2$ as it passes through face-centres or along edges. (The map $\{4,4\}_{b,c}$, for integers $b, c\ge 0$, is the quotient of the unit square tessellation of $\mathbb C$ by the translation group corresponding to the ideal $(b+ci)$ in the ring ${\mathbb Z}[i]$ of Gaussian integers; see~\cite[\S8.3]{CM} for its properties.) For brevity we will simply write that a reflection $r$ has type $i$ if it has an action of type $i$ on some flag, noting that a reflection may have more than one type.

\section{Conjugacy classes of reflections of maps}

If $H$ is a group of automorphisms of a map $\mathcal M$, we will denote the number of conjugacy classes of reflections of $\mathcal M$ in $H$ by $cr(H)$, and the number of classes of reflections of type $i=0, 1$ or $2$ by $cr_i(H)$. If $H\le\Gamma$ it will be assumed (unless stated otherwise) that ${\mathcal M}={\mathcal M}_{\infty}$. We will also use the notations $cr({\mathcal M})$ and $cr_i({\mathcal M})$ in cases where $H={\rm Aut}\,{\mathcal M}$ for some map $\mathcal M$.

If $\mathcal M$ is the map corresponding a a subgroup $M\le\Gamma$, than any reflection $r$ of $\mathcal M$ is induced by a reflection $R\in N:=N_{\Gamma}(M)$. This must be a conjugate $R_i^C$, where $C\in\Gamma$, of one of the standard generators $R_i\;(i=0, 1, 2)$ of $\Gamma$. We now consider how such reflections are divided into conjugacy classes in $N$.

\begin{prop}\label{boundcrN} Let $N$ be any subgroup of $\Gamma$. Then:
\begin{enumerate}

\item $cr_0(N)$ is the number of cycles of $R_2$ on those cosets of $N$ in $\Gamma$ fixed by $R_0$; each conjugacy class of type~0 reflections in $N$ consists of the elements $CR_0C^{-1}$ where $C$ lies in a particular coset of $N$ fixed by $R_0$ and $R_2$, or in a particular pair of cosets fixed by $R_0$ and transposed by $R_2$.

\item $cr_1(N)$ is the number of cosets of $N$ in $\Gamma$ fixed by $R_1$; each conjugacy class of type~1 reflections in $N$ consists of the elements $CR_1C^{-1}$ where $C$ lies in a particular coset fixed by $R_1$.

\item $cr_2(N)$ is the number of cycles of $R_0$ on those cosets of $N$ in $\Gamma$ fixed by $R_2$; each conjugacy class of type~2 reflections in $N$ consists of the elements $CR_2C^{-1}$ where $C$ lies in a particular coset of $N$ fixed by $R_0$ and $R_2$, or in a particular pair of cosets fixed by $R_2$ and transposed by $R_0$.

\end{enumerate} 
\end{prop}

(Note that the cosets of $N$ in $\Gamma$ are in bijective correspondence with the flags of the map ${\mathcal N}={\mathcal M}/{\rm Aut}\,{\mathcal M}$ with map subgroup $N$, permuted by the generators $R_i$ of $\Gamma$ in the standard way, so the numbers $cr_i(N)$ can be found by considering how vertices, edges and flags of $\mathcal N$ meet the boundary of $\mathcal N$.)

\medskip

\noindent{\sl Proof.}  One can write $\Gamma$ as a disjoint union of cosets $NC_j$ of $N$, where $j$ ranges over some index set $J\;(=\{1,\ldots, n\}$ if $|\Gamma:N|=n$ is finite), with $N$ represented by $C_1=1$. The action of $\Gamma$ on the cosets of $N$ can then be identified with its action on $J$, so that $N$ is the subgroup of $\Gamma$ fixing $1$.

Any reflection $R\in \Gamma$ must be have the form $R_i^C$ for some $i=0, 1$ or $2$ and $C\in\Gamma$.  Then $R\in N$ if and only if $R$ fixes $1\in J$, or equivalently $R_i$ fixes the image $j\in J$ of $1$ under $C^{-1}$. The set of all elements sending $1$ to $j$ is the coset $NC_j$, so $C^{-1}\in NC_j$ and thus $R$ is conjugate in $N$ to $C_jR_iC_j^{-1}$. Thus the conjugacy classes in $N$ of reflections of type $i$ are represented by elements $C_jR_iC_j^{-1}$, where $R_i$ fixes $j\in J$.

Since $R_i$ is not conjugate in $\Gamma$ to $R_{i'}$ for $i\ne i'$, it follows that  $C_jR_iC_j^{-1}$ cannot be conjugate in $N$ to $C_kR_{i'}C_k^{-1}$ unless $i=i'$ and $R_i$ fixes both $j$ and $k$. Such a conjugacy is equivalent to some element of $C_k^{-1}NC_j$ centralising $R_i$, that is, to some element of the centraliser $C_{\Gamma}(R_i)$ of $R_i$ in $\Gamma$ sending $k$ to $j$. 

Suppose first that $i=1$. The normal form theorem for free productsm (see~\cite[\S IV.1]{LS} or~\cite[\S 4.1]{MKS}) implies that $C_{\Gamma}(R_1)=\langle R_1\rangle$, and since we are assuming that $R_1$ fixes $j$ and $k$ we have $j=k$. Thus the conjugacy classes of reflections of type~1 in $N$ are in bijective correspondence with the fixed points of $R_1$ on $J$, and hence on the cosets of $N$.

Now suppose that $i=0$. In this case $C_{\Gamma}(R_0)=\langle R_0, R_2\rangle$, with $R_0$ fixing $j$ and $k$, so either $j=k$ or $R_2$ (and $R_0R_2$) transpose $j$ and $k$. Thus the conjugacy classes of reflections of type~0 in $N$ are in bijective correspondence with the cycles of $R_2$ on the fixed points of $R_0$. By applying the automorphism of $\Gamma$ transposing $R_0$ and $R_2$, or equivalently by map duality, we obtain the corresponding result for reflections of type~2. \hfill$\square$

\begin{cor}\label{boundcrM}
Let $\mathcal M$ be a map with map subgroup $M\le\Gamma$, and let $N=N_{\Gamma}(M)$. Then
\begin{enumerate}
\item $cr_0({\mathcal M})$ is less than or equal to the number of cycles of $R_2$ on the fixed points of $R_0$ on the cosets of $N$ in $\Gamma$;
\item $cr_1({\mathcal M})$ is less than or equal to the number of fixed points of $R_1$ on the cosets of $N$ in $\Gamma$;
\item $cr_2({\mathcal M})$ is less than or equal to the number of cycles of $R_0$ on the fixed points of $R_2$ on the cosets of $N$ in $\Gamma$.
\end{enumerate}
\end{cor}

\noindent{\sl Proof.} The reflections of type~$i$ in ${\rm Aut}\,{\mathcal M}\cong N/M$ are the images of those in $N$, and conjugate reflections in $N$ have conjugate images, so
\begin{equation}\label{crMcrN}
cr_i({\mathcal M})\le cr_i(N).
\end{equation}
This inequality, together with Proposition~\ref{boundcrN}, gives the required upper bound on $cr_i({\mathcal M})$.\hfill$\square$

\medskip

The inequality in~(\ref{crMcrN}) may be strict in some cases, since nonconjugate reflections in $N$ could have conjugate images in ${\rm Aut}\,{\mathcal M}$. Nevertheless, each reflection of type~$i$ in $N$ induces a reflection of the same type in ${\rm Aut}\,{\mathcal M}$, so we have

\begin{lemma}
Let $\mathcal M$ and $N$ be as above. If $cr_i(N)\ge 1$ then $cr_i({\mathcal M})\ge 1$.
\end{lemma}

Similarly, although
\begin{equation}
cr(N)=\sum_{i=0}^2cr_i(N),
\end{equation}
the inequality
\begin{equation}\label{crMbd}
cr({\mathcal M})\le \sum_{i=0}^2cr_i({\mathcal M})
\end{equation}
could be strict, since a reflection of $\mathcal M$ could have more than one type whereas each reflection in $N$ has a unique type.

\medskip

\noindent{\bf Example 4.1} If $\mathcal M$ is regular then $N=\Gamma$, so $cr_i({\mathcal M})=cr_i(\Gamma)=1$ for each $i$ and hence $cr({\mathcal M})\le 3$. On the other hand, if $\mathcal M$ is orientably regular but not regular then $N=\Gamma^+$, so $cr_i({\mathcal M})=cr_i(\Gamma^+)=0$ for each $i$ and hence $cr({\mathcal M})=0$.

\medskip

Note that if $cr_1(N)$ is odd then $|\Gamma:N|$ is odd, so $R_0$ and $R_2$ both have fixed points, and hence $cr_i(N)\ge 1$ for $i=0, 2$. In fact, this is the only restriction on the values of the integers $cr_i(N)$:

\begin{thm}\label{criN}
Given any integers $c_0, c_1, c_2\ge 0$, with $c_0, c_2\ge 1$ if $c_1$ is odd, there is a subgroup $N$ of finite index in $\Gamma$ with $cr_i(N)=c_i$ for $i=0, 1$ and $2$.
\end{thm}

\noindent{\sl Proof.} We will take $N$ to be the stabiliser of a point in a suitable permutation representation $\theta:\Gamma\to S_n$ of $\Gamma$. We will construct $\theta$ by first specifying the action of $E=\langle R_0, R_2\rangle$, and then that of $R_1$.

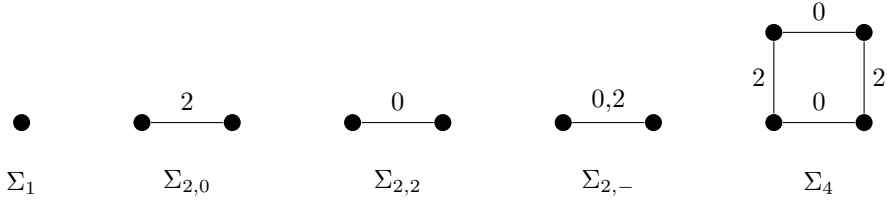
\begin{figure}[h!]
\label{K5}
\begin{center}
 \begin{tikzpicture}[scale=0.4, inner sep=0.8mm]
 
 \node at (-12,0) [shape=circle, fill=black] {};
 \node at (-12,-2) {$\Sigma_1$};
 

\node (a) at (-8,0) [shape=circle, fill=black] {};
\node (b) at (-5,0) [shape=circle, fill=black] {};
\draw (a) to (b);
\node at (-6.5,0.7) {2};
 \node at (-6.5,-2) {$\Sigma_{2,0}$};


\node (c) at (-1,0) [shape=circle, fill=black] {};
\node (d) at (2,0) [shape=circle, fill=black] {};
\draw (c) to (d);
\node at (0.5,0.7) {0};
 \node at (0.5,-2) {$\Sigma_{2,2}$};


\node (e) at (6,0) [shape=circle, fill=black] {};
\node (f) at (9,0) [shape=circle, fill=black] {};
\draw (e) to (f);
\node at (7.5,0.7) {0,2};
 \node at (7.5,-2) {$\Sigma_{2,-}$};


\node (g) at (13,0) [shape=circle, fill=black] {};
\node (h) at (16,0) [shape=circle, fill=black] {};
\draw (g) to (h);
\node (i) at (13,3) [shape=circle, fill=black] {};
\node (j) at (16,3) [shape=circle, fill=black] {};
\draw (i) to (j);
\draw (g) to (i);
\draw (h) to (j);
\node at (14.5,0.7) {0};
\node at (14.5,3.7) {0};
\node at (12.5,1.5) {2};
\node at (16.5,1.5) {2};
 \node at (14.5,-2) {$\Sigma_4$};

\end{tikzpicture}

\end{center}
\caption{Permutation graphs for orbits of $E=\langle R_0, R_1\rangle$} 
\end{figure}

Any orbit of $E$ has length $l=1, 2$ or $4$, and will accordingly be denoted by $\Omega_l$. If $l=1$ or $4$ then there is, up to isomorphism, a unique action of $E$ on $\Omega_l$; we will denote an orbit of length $l=2$ by $\Omega_{2,0}$, $\Omega_{2,2}$ or $\Omega_{2,-}$ as the kernel of the action of $E$ is generated by $R_0$, $R_2$ or $R_0R_2$ respectively. Any orbit $\Omega_1$ provides a single fixed point for each of $R_0$ and $R_2$, giving two reflections conjugate to $R_0$ and $R_2$, while any orbit $\Omega_{2,i}\;(i=0, 2)$ provides a $2$-cycle for $R_{2-i}$ on the fixed points of $R_i$, and hence a reflection conjugate to $R_i$; other orbits of $E$ provide no fixed points for these two generators, and hence no reflections. Let $\Sigma_1$, $\Sigma_{2,i}$, $\Sigma_{2,-}$ and $\Sigma_4$ denote the permutation graphs for these orbits of $E$ with respect to the generators $R_0, R_2$, with their $2$-cycles represented as undirected edges labelled $0$, $2$ or both, but omitting loops for fixed points (see Figure~7). We will construct a permutation graph $\Sigma$ for $\Gamma$ as a union of such subgraphs, connected by edges labelled $1$ representing $2$-cycles of $R_1$.

First suppose that $c_1$ is even. Take $c_0$ copies of $\Sigma_{2,0}$, together with $c_1/2$ copies of $\Sigma_4$, and $c_2$ copies of $\Sigma_{2,2}$. Now arrange these subgraphs in any cyclic order, and join each neighbouring pair with a single undirected edge labelled $1$; in the case of subgraphs $\Sigma_4$ we may attach these edges to any pair of vertices, leaving the other two as fixed points of $R_1$. The resulting `necklace' graph $\Sigma$ is connected, so it is the permutation graph for a transitive representation $\theta:\Gamma\to S_n$ of $\Gamma$ where $n=2(c_0+c_1+c_2)$, or equivalently the Schreier coset graph for a point stabiliser $N$, a subgroup of index $n$ in $\Gamma$. For $i=0, 2$ the $c_i$ subgraphs $\Sigma_{2,i}$ provide $c_i$ cycles of $R_{2-i}$ on the fixed points of $R_i$, and the $c_1/2$ subgraphs $\Sigma_4$ provide $c_1$ fixed points of $R_1$, two from each subgraph. The generators $R_i$ have no other fixed points, so Proposition~\ref{boundcrN} implies that $cr_i(N)=c_i$ for each $i=0, 1, 2$.

If $c_1$ is odd, then $c_0, c_2\ge 1$ by our hypothesis. In this case take $c_0-1$, $(c_1+1)/2$ and $c_2-1$ subgraphs $\Sigma_{2,0}$, $\Sigma_4$ and $\Sigma_{2,2}$, then use edges labelled $1$ to join these in any cyclic order as above, and finally to join a single copy of $\Sigma_1$ to a vertex in one of the subgraphs $\Sigma_4$ not yet incident with such an edge. As before the resulting graph $\Sigma$ defines a transitive permutation representation of $\Gamma$, this time of degree $n=2(c_0+c_1+c_2)-1$, such that a point stabiliser $N$ satisfies $cr_i(N)=c_i$ for each $i$.\hfill$\square$

\medskip

\noindent{\bf Remark} Note, for future use, that in either part of this proof we could also have included in the necklace any number of subgraphs $\Sigma_{2,-}$, and any number of subgraphs $\Sigma_4$ each containing a single edge labelled $1$, since these contribute no further fixed points of any $R_i$.

\medskip

In order to prove an analogous result for maps, we need the following concept. A group $H$ is {\em conjugacy separable\/} if, whenever two elements are not conjugate in $H$, there is some finite quotient of $H$ in which their images are also not conjugate. This property extends to all finite sets of mutually non-conjugate elements. Finite groups and finitely generated abelian groups are conjugacy separable, as are free products of conjugacy separable groups~\cite{Rem, Ste}. This immediately implies that $\Gamma$ and its subgroups $N$ are conjugacy separable, so $N$ has finite quotients $N/M$ such that distinct conjugacy classes of reflections in $N$ have distinct images in $N/M$. Here $N_{\Gamma}(M)\ge N$, but in order to apply this to maps we need to ensure that $N_{\Gamma}(M)=N$; for this we need more explicit constructions of $N$ and $M$.

\begin{cor}
Given any integers $c_0, c_1, c_2\ge 0$, with $c_0, c_2\ge 1$ if $c_1$ is odd, there is a finite map $\mathcal M$ with $cr_i({\mathcal M})=c_i$ for $i=0, 1$ and $2$, and $cr({\mathcal M})=c_0+c_1+c_2$.
\end{cor}

\noindent{\sl Proof.} Let $N$ be the subgroup of finite index $n$ in $\Gamma$ constructed in the proof of Theorem~\ref{criN}, with $cr_i(N)=c_i$ for each $i$, so that $cr(N)=c_0+c_1+c_2$. As noted earlier, since $\Gamma=E*\langle R_1\rangle$ the Kurosh Subgroup Theorem implies that $N$ is a free product of a free group and subgroups which are conjugate to $\langle R_1\rangle$ or to subgroups $F\le E$, corresponding to the flags of $\mathcal N$ with these finite groups as stabilisers. Such flags are as described in the proof of Theorem~\ref{criN}.

Let $M=N'N^2$, the group generated by the commutators and squares of the elements of $N$. This is normal subgroup of finite index in $N$, and $N/M$ is an elementary abelian $2$-group. Non-conjugate reflections in $N$ will have distinct (and thus non-conjugate) involutions as their images in $N/M$. By construction, $N\le N_{\Gamma}(M)$, and if we can prove equality here then the map $\mathcal M$ corresponding to $M$ will have the required properties.

First suppose that $c_1$ is odd, so that $c_0, c_2\ge 1$ and the index $n=|\Gamma:N|$ is odd. We can modify the construction of the permutation representation $\theta:\Gamma\to S_n$ in the proof of Theorem~\ref{criN} by including sufficiently many copies of the subgraph $\Sigma_{2,-}$ in the necklace, so that the degree $n$ is a prime $p$. This does not change the values of $cr_i(N)$ or $cr(N)$, and it ensures that $\theta$ is a primitive representation. The point stabiliser $N$ is therefore a maximal subgroup of $\Gamma$, so $N_{\Gamma}(M)$ is either $N$ or $\Gamma$. In the latter case, as a normal subgroup of $\Gamma$ contained in $N$, $M$ must be contained in the core $K={\rm ker}\,\theta$ of $N$ in $\Gamma$. It follows that $N/K$ is also an elementary abelian $2$-group. Now $N/K$ is isomorphic to the stabiliser $G_0$ of a point in the permutation group $G={\rm Mon}\,{\mathcal N}\cong\Gamma/K$ of degree $n=p$ induced by $\Gamma$. Since $G$ has order divisible by only two primes ($2$ and $p$), Burnside's $p^aq^b$ theorem (see~\cite{Bur} or~\cite[Hauptsatz~V.7.3]{Hup}) implies that it is solvable. As a solvable group of prime degree $p$, $G$ is isomorphic to a subgroup of $AGL_1(p)$ by a theorem of Galois (see~\cite[Satz~II.3.6]{Hup}), so $G_0$ is cyclic. This is impossible, since the Klein four-group $E$ acts faithfully and fixes a point (the single vertex in the subgraph $\Sigma_1$ of $\Sigma$). Thus $N_{\Gamma}(M)=N$, so the map $\mathcal M$ corresponding to $M$ has its conjugacy classes of reflections corresponding to those in $N$, as required. 

Now suppose that $c_1$ is even, so the index $n=|\Gamma:N|$ is also even. Let $\Sigma_4^*$ denote a copy of $\Sigma_4$ with an extra edge labelled $1$ joining two non-adjacent vertices. As before, we can modify the construction of $\theta$ by adding copies of $\Sigma_{2,-}$, and now also of $\Sigma_4^*$, to the necklace, this time ensuring that $n=2p$ for some odd prime $p$. Again, let $M=N'N^2$. If $N_{\Gamma}(M)=N$ we are done, so we may assume that $N_{\Gamma}(M)>N$. If $\theta$ is primitive then as before we find that $\Gamma$ induces a solvable subgroup of $S_n$, whereas primitive solvable groups have prime power degree. We may therefore assume that $\theta$ is imprimitive, so there are either $p$ blocks of size $2$, or $2$ of size $p$, corresponding to a subgroup $L$ of $\Gamma$ containing $N$ with index $|L:N|=2$ or $p$. In the first case, $N$ is normal in $L$, so $L$ induces an automorphisms of order $2$ on the map $\mathcal N$ corresponding to $N$, and hence on the labelled graph $\Sigma$; however, we can arrange the cyclic order of subgraphs around the necklace so that no such automorphism exists. In the second case, $L$ has index $2$ in $\Gamma$, and contains all elements of $\Gamma$ with fixed points under $\theta$. Including subgraphs $\Sigma_{2,-}$ and $\Sigma_4^*$ in the necklace ensures that $R_0R_2\in L$ and $R_1R_0R_2\in L$, so $R_1\in L$. Since $R_0, R_1$ and $R_2$ generate $\Gamma$, this implies that $R_0, R_2\not\in L$, so $c_0=c_2=0$. Thus the subgraphs in the necklace all have the form $\Sigma_{2,-}$, $\Sigma_4^*$ or $\Sigma_4$. There are $c_1/2$ of the latter, so if there are $a$ and $b$ of the former, then $n=2a+4b+2c_1$, giving $a+2b+c_1=p$, so $a$ is odd. However, there is a closed walk around the necklace corresponding to a word $w(R_0, R_1, R_2)\in N$ in which $a+b+c_1/2$ instances of $R_1$ alternate with a total of $a$ instances of $R_0$ and $b+c_1/2$ of $R_0R_2$. Since $R_1, R_0R_2\in L$ this gives $R_0^a\in L$ and hence $R_0\in L$, a contradiction.\hfill$\square$

\medskip

Note that, unlike the situation with $cr_i(N)$, if $cr_1({\mathcal M})$ is odd for some map $\mathcal M$ it does not follow that $cr_i({\mathcal M})\ge 1$ for $i=0, 2$: the point is that $cr_1(N)$ could be greater than $cr_1({\mathcal M})$, and in particular could be even, allowing $cr_i(N)=0$ and hence $cr({\mathcal M})=0$ for $i=0$ or $2$. We will construct examples in Section~5.6.

\section{Reflections of edge-transitive maps}

In contrast with the earlier results for vertex- and face-transitive maps, where there can be any number of conjugacy classes of reflections, the situation for edge-transitive maps is quite different. The main aim of this section is to prove Theorem~\ref{edgetrans}, giving the upper bound $cr({\mathcal M})\le 4$ for such maps and describing those maps attaining it. 

\subsection{Types of reflections of edge-transitive maps}

At first glance it might appears that an edge-transitive map $\mathcal M$ could have as many as six conjugacy classes of reflections, acting on the flags at a particular edge $e$ as follows:

\begin{itemize}
\item[a)] one each of type $0$ or $2$, reversing or fixing $e$, and
\item[b)] four of type $1$, transposing $e$ with one of the four edges sharing a common vertex and face with $e$.
\end{itemize}
However, as noted earlier, a reflection may have actions of different types on different flags: for instance, if $\mathcal M$ has a vertex or face of odd valency then by edge-transitivity one of these is incident with $e$, so any reflection of type $2$ or $0$ respectively at $e$ will have an action of type $1$ at some other edge; by edge-transitivity, this implies that $r$ is conjugate to a reflection of type $1$ at $e$, thus reducing the number of conjugacy classes. Even if this problem is avoided, for instance by considering maps with only even valencies, there is a further problem, this time unavoidable: if $e$ admits a reflection of type~$0$ or $2$ then this will conjugate any reflection of type~$1$ at $e$ to another of type $1$, again reducing the number of classes of reflections. In fact, by using this argument it is easy to see that $cr(\mathcal M)\le 4$, and that this bound is attained if and only if, for some (and hence every) edge $e$, all four reflections of type~$1$ are automorphisms of $\mathcal M$, and no two of them are conjugate to each other. This last condition implies that $\mathcal M$ has no reflections of type $0$ or $2$ as automorphisms, so that $\mathcal M$ has automorphism type~3 in the taxonomy of edge-transitive maps described by Graver and Watkins~\cite[Table~2]{GW}. It also implies that all vertex- and face-valencies are divisible by $4$.

\subsection{Taxonomy of edge-transitive maps}

A map $\mathcal M$ with map subgroup $M\le\Gamma$ is edge-transitive if and only if $\Gamma=NE$, where $N:=N_{\Gamma}(M)$ and $E:=\langle R_0, R_2\rangle\cong V_4$ is the stabiliser in $\Gamma$ of an edge. This is equivalent to $E$ acting transitively on the cosets of $N$, so the index $n:=|\Gamma:N|$ of $N$ in $\Gamma$ divides $|E|=4$, giving $n=1, 2$ or $4$.

It is a simple matter to find all such subgroups $N$ of $\Gamma$: they are the stabilisers of points in those permutation representations $\Gamma\to S_n$ of $\Gamma$ such that $E$ acts transitively. These representations can be determined by considering the possible images of the generators $R_i$ of $\Gamma$ in $S_n$. Case-by-case analysis shows that there are the following $1+6+7=14$ conjugacy classes of such subgroups $N$, listed below with the automorphism types of the corresponding maps $\mathcal M$ as assigned by Graver and Watkins in~\cite[Table~2]{GW}; the presence or absence of the symbols $V, F, P$ indicates whether or not ${\rm Aut}\,{\mathcal M}$ is transitive on vertices, faces or Petrie polygons (closed zig-zag paths, turning right and left at alternate vertices~\cite[\S5.2, \S8.6]{CM}).
In each case one can use Proposition~\ref{boundcrN} to determine $cr_i(N)$ for each $i=0, 1, 2$, and hence to obtain information about the possible values of $cr_i({\mathcal M})$ for maps $\mathcal M$ with $N_{\Gamma}(M)=N$.

\begin{itemize}

\item {\bf Type $1$} ($V, F, P$). If $n=1$ then $N=\Gamma$, corresponding to $\mathcal M$ being regular, with $cr_i(N)=1$ for $i=0, 1, 2$.  It follows that $cr_i({\mathcal M})=1$ for each $i$, so $1\le cr({\mathcal M})\le 3$;  the reflections of $\mathcal M$ are induced by the conjugates of  $R_0, R_1$ and $R_2$ in $\Gamma$.

\end{itemize}

If $n=2$ there are six possible subgroups $N$, each of them normal in $\Gamma$, namely the seven subgroups of index $2$ except $\Gamma'E$. The different cases are distinguished by which generators $R_i$ are in $N$, as follows:

\begin{itemize}

\item {\bf Type $2^P$ex} ($V, F, P$). No $R_i\in N$, that is, $N=\Gamma^+=\langle R_0R_1, R_1R_2\rangle$ or equivalently $\mathcal M$ is orientably regular and chiral; then $cr_i({\mathcal M})=cr_i(N)=0$ for each $i$, and there are no reflections.

\item {\bf Type $2^*$ex} ($V, F, P$). Only $R_0\in N$; $cr_i({\mathcal M})=cr_i(N)=0$ for $i=1, 2$ and $cr_0({\mathcal M})=cr_0(N)=1$, so $cr({\mathcal M})=cr(N)=1$; the reflections are induced by the conjugates of $R_0$ in $N$. (Note that $R_2R_1\in N$, so $R_0^{R_1}=R_0^{R_2R_1}$ is conjugate in $N$ to $R_0$.)

\item  {\bf Type $2^P$} ($V, F$). Only $R_1\in N$; $cr_i({\mathcal M})=cr_i(N)=0$ for $i=0, 2$ and $1\le cr_1({\mathcal M})\le cr_1(N)=2$, so $1\le cr({\mathcal M})\le cr(N)=2$; the reflections are induced by the conjugates of $R_1$ and $R_1^{R_0}$ in $N$. (Note that $R_0R_2\in N$, so $R_1^{R_2}$ and $R_1^{R_0R_2}$ are conjugate to $R_1^{R_0}$ and $R_1$ respectively.)

\item  {\bf Type $2$ex} ($V, F, P$). Only $R_2\in N$;  $cr_i({\mathcal M})=cr_i(N)=0$ for $i=0, 1$ and $cr_2({\mathcal M})=cr_2(N)=1$, so $cr({\mathcal M})=cr(N)=1$; the reflections are induced by the conjugates of $R_2$ in $N$.

\item  {\bf Type $2^*$} ($V, P$). Only $R_0, R_1\in N$;  $cr_0({\mathcal M})=cr_0(N)=1$, $1\le cr_1({\mathcal M})\le cr_1(N)=2$ and $cr_2({\mathcal M})=cr_2(N)=0$, so $1\le cr({\mathcal M})\le cr(N)=3$; the reflections are induced by the conjugates of $R_0, R_1$ and $R_1^{R_2}$ in $N$.

\item  {\bf Type $2$} ($F, P$). Only $R_1, R_2\in N$; $cr_0({\mathcal M})=cr_0(N)=0$,  $1\le cr_1({\mathcal M})\le cr_1(N)=2$ and $cr_2({\mathcal M})=cr_2(N)=1$, so $1\le cr({\mathcal M})\le cr(N)=3$; the reflections are induced by the conjugates of $R_1, R_1^{R_0}$ and $R_2$ in $N$.

\end{itemize}

If $n=4$ there are seven conjugacy classes of subgroups $N$, corresponding to permutation representations $\Gamma\to S_4$ with $E$ acting regularly. Without loss of generality we may suppose that $R_0\mapsto(12)(34)$ and $R_2\mapsto(14)(23)$, so (up to conjugation fixing $R_0$ and $R_2$, that is, by elements of $E$) the possibilities for the image of $R_1$ in $S_4$ are as follows:

\begin{itemize}

\item  {\bf Type $3$} ($-$). $R_1\mapsto$ the identity, so $\Gamma$ acts as $V_4$ on the cosets of $N$, which is normal in $\Gamma$; then $cr_i({\mathcal M})=cr_i(N)=0$ for $i=0, 2$ while $1\le cr_1({\mathcal M})\le cr_1(N)=4$, so $1\le cr({\mathcal M})\le cr(N)=4$; the conjugacy classes of reflections are represented by $R_1, R_1^{R_0}, R_1^{R_2}$ and $R_1^{R_0R_2}$. 

\item {\bf Types $4$} ($F, P$), {\bf $4^P$} ($V, F$) and {\bf $4^*$} ($V, P$). $R_1\mapsto(23), (24)$ or $(34)$, so $\Gamma$ acts as $D_4$ on the cosets of $N$, which has two conjugates in $\Gamma$; then $cr_i({\mathcal M})=cr_i(N)=0$ for $i=0, 2$ while $1\le cr_1({\mathcal M})\le cr_1(N)=2$, so $1\le cr({\mathcal M})\le cr(N)=2$; the conjugacy classes of reflections are represented by $R_1$ and either $R_1^{R_0}$, $R_1^{R_0R_2}$ or $R_1^{R_2}$ respectively.

\item  {\bf Types $5^*$} ($V, P$), {\bf $5^P$} ($V, F$) and {\bf $5$} ($F, P$). $R_1\mapsto(12)(34), (13)(24)$ or $(14)(23)$, so $\Gamma$ acts as $V_4$ on the cosets of $N$, which is normal in $\Gamma$; then $cr_i({\mathcal M})=cr_i(N)=0$ for each $i$, and there are no reflections.

\end{itemize}

\subsection{Upper bounds}

Inspection of these $14$ cases proves the upper bound $cr({\mathcal M})\le 4$ for edge-transitive maps, stated in Theorem~\ref{edgetrans}, and shows that, as claimed there, it is attained only by certain maps of type~3. In each of the $14$ cases one can easily draw the map $\mathcal N$ corresponding to $N$; the edge-transitive maps $\mathcal M$ of that type are all regular coverings of this map, with covering group ${\rm Aut}\,{\mathcal M}\cong N/M$. For instance, the maps of type~3, considered in the next section, are all regular coverings of the map on the disc shown in Figure~5.

Note also that the orientably regular maps, those with $N\ge\Gamma^+$, are those of type~1 (the regular maps) or~$2^P$ex (the chiral maps), so this also confirms the upper bound $cr({\mathcal M})\le 3$ for such maps stated in the first sentence of this paper. For instance, the tetrahedron and the cube, regarded as orientably regular maps on the sphere, have $cr({\mathcal M})=1$ and $2$, while the torus maps $\{4,4\}_{2,0}$ and $\{4,4\}_{2,1}$ have $cr({\mathcal M})=3$ and $0$. These are shown in Figure~8, where in each case opposite sides of the outer square are identified to form a torus; in the first case, the axes of three non-conjugate reflections are indicated by broken lines.

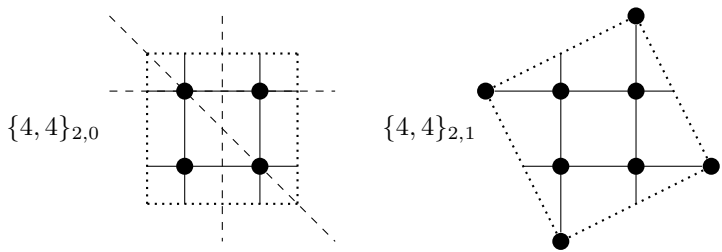
\begin{figure}[h!]
\label{K5}
\begin{center}
 \begin{tikzpicture}[scale=0.5, inner sep=0.8mm]

\node (A) at (-6,6) [shape=circle, fill=black] {};
\node (B) at (-4,6) [shape=circle, fill=black] {};
\node (C) at (-4,4) [shape=circle, fill=black] {};
\node (D) at (-6,4) [shape=circle, fill=black] {};

\draw (-7,6) to (-3,6);
\draw (-4,7) to (-4,3);
\draw (-3,4) to (-7,4);
\draw (-6,3) to (-6,7);

\draw [dotted, thick] (-7,7) to (-3,7);
\draw [dotted, thick] (-3,7) to (-3,3);
\draw [dotted, thick] (-3,3) to (-7,3);
\draw [dotted, thick] (-7,3) to (-7,7);

\draw [dashed] (-5,8) to (-5,2);
\draw [dashed] (-8,8) to (-2,2);
\draw [dashed] (-8,6) to (-2,6);

\node at (-9.5,5) {$\{4,4\}_{2,0}$};

\node (a) at (6,8) [shape=circle, fill=black] {};
\node (b) at (6,6) [shape=circle, fill=black] {};
\node (c) at (4,6) [shape=circle, fill=black] {};
\node (d) at (2,6) [shape=circle, fill=black] {};
\node (e) at (8,4) [shape=circle, fill=black] {};
\node (f) at (6,4) [shape=circle, fill=black] {};
\node (g) at (4,4) [shape=circle, fill=black] {};
\node (h) at (4,2) [shape=circle, fill=black] {};

\draw (7,6) to (d);
\draw (e) to (3,4);
\draw (a) to (6,3);
\draw ((4,7) to (h);

\draw [dotted, thick] (a) to (d);
\draw [dotted, thick] (d) to (h);
\draw [dotted, thick] (h) to (e);
\draw [dotted, thick] (e) to (a);

\node at (0.5,5) {$\{4,4\}_{2,1}$};

\end{tikzpicture}

\end{center}
\caption{The torus maps $\{4,4\}_{2,0}$ and $\{4,4\}_{2,1}$} 
\end{figure}

\subsection{Examples attaining the upper bound}

The following family of examples shows that the upper bound $cr({\mathcal M})\le 4$ in Theorem~\ref{edgetrans} is attained.

\medskip

\noindent{\bf Example 5.1} Let $\mathcal T$ be the torus map $\{4,4\}_{b,0}$ for some integer $b\ge 1$ (see Figure~8 for the case $b=2$). This is a regular map with $b^2$ vertices and faces, $2b^2$ edges, and $8b^2$ automorphisms~\cite[\S 8.3]{CM}. If $b$ is even then $\mathcal T$ is bipartite and $2$-face-colourable, meaning that its vertices and faces can be $2$-coloured so that any pair of vertices or pair of faces sharing an edge have different colours.

Now let $\mathcal M$ be a double covering of $\mathcal T$, branched over alternate vertices and alternate faces, in each case those of a particular colour. This map has $3b^2/2$ vertices and faces, and $4b^2$ edges, so it has Euler characteristic $-b^2$ and hence its genus is $1+b^2/2$. Let $H$ be the subgroup of index $4$ in ${\rm Aut}\,{\mathcal T}$ preserving the colours of the vertices and faces of $\mathcal T$. This group, generated by the four reflections of type $1$ associated with a given edge of $\mathcal T$, acts regularly on the edges of $\mathcal T$, and it lifts to a group $G\le{\rm Aut}\,{\mathcal M}$ acting regularly on the edges of $\mathcal M$, so $\mathcal M$ is edge-transitive. The two vertices and the two faces of $\mathcal M$ incident with a given edge $e$ of $\mathcal M$ have different valencies (namely $8$ and $4$), so $e$ has trivial stabiliser in ${\rm Aut}\,{\mathcal M}$ and hence ${\rm Aut}\,{\mathcal M}=G$. None of the four type $1$ reflections of $\mathcal M$ at $e$ can be conjugate to each other in $G$: if they were, their images in $H$ would be conjugate in $H$, whereas these lie in distinct conjugacy classes in $H$, distinguished by the colours of the vertices and faces they leave invariant (since their fixed-point sets in $\mathcal T$ are connected). Thus $cr({\mathcal M})=4$.

\subsection{Just-edge-transitive maps}

In Section~5.2, the only case allowing the possibility that $cr({\mathcal M})=4$ is that in which the maps $\mathcal M$ have automorphism type~3, so that $N$ is the normal closure $N_3:=\langle\langle R_1\rangle\rangle^{\Gamma}$ of $R_1$ in $\Gamma$. These are the just-edge-transitive maps, those that are edge- but neither vertex- nor face-transitive, studied in~\cite{Jon}. Equivalently, the map $\mathcal N$ corresponding to $N$ is the map ${\mathcal N}_3$ on the closed disc, with two vertices, one edge and two faces, shown in Figure~5 and described in Example~3.1. As shown there, this subgroup $N_3$ has generators
\[S_1=R_1,\; S_2=R_1^{R_2},\; S_3=R_1^{R_0R_2}=R_1^{R_2R_0}\quad{\rm and}\quad S_4=R_1^{R_0},\]
with defining relations $S_i^2=1$ for $i=1,\ldots, 4$, so $N_3$ is the free product of four copies of $C_2$. The corresponding quotients $G=N_3/M\cong{\rm Aut}\,{\mathcal M}$ are those groups which can be generated by four elements $s_i$ (the images of $S_i$) satisfying $s_i^2=1$. Now
\[S_1S_2=R_1R_1^{R_2}=(R_1R_2)^2,\]
and
\[S_3S_4=R_1^{R_2R_0}R_1^{R_0}=((R_1R_2)^{R_0})^2,\]
so if each $s_is_j$ has order $p_{ij}\;(=p_{ji})$ then the two vertices incident with each edge have valencies $2p_{12}$ and $2p_{34}$; similarly the incident faces are a $2p_{14}$-gon and a $2p_{23}$-gon, while the incident Petrie polygons have lengths $2p_{13}$ and $2p_{24}$. This is equivalent to factoring the epimorphism $N_3\to G$ through the rank $4$ Coxeter group $C=C(p_{ij})$ obtained by adding the relations $(S_iS_j)^{p_{ij}}=1$ to the presentation of $N_3$, and requiring the epimorphism $C\to G$ to map the dihedral subgroups $\langle S_i, S_j\rangle$ of $C$ faithfully into $G$. The existence of such an epimorphism onto a finite group $G$ follows from the residual finiteness of Coxeter groups, as finitely generated linear groups (by Mal'cev's Theorem~\cite{Mal}). Reflections $s_i$ and $s_j$ of $\mathcal M$ are then conjugate in $G$ if and only if $S_i$ and $S_j$ are conjugate in $C$, that is, if and only if $i$ and $j$ are equivalent in the transitive closure of the relation $\{(i, j)\mid p_{ij}\;\hbox{is odd}\}$, so we have $cr({\mathcal M})=cr(C)$.

For any normal subgroup $M$ of $N_3$ we have $N_{\Gamma}(M)\ge N_3$, and for $\mathcal M$ to be just-edge-transitive we need $N_{\Gamma}(M)=N_3$. Since each of the representatives $R_0, R_2$ and $R_0R_2$ of the nontrivial cosets of $N_3$ in $\Gamma$ acts by conjugation as a double transposition on the generators $S_i$ of $N_3$, this is equivalent to the property that no automorphism of $G$ induces a double transposition of its generators $s_i$. A simple way of achieving this is to ensure that at least two of the conditions $p_{12}\ne p_{34}$, $p_{13}\ne p_{24}$ and $p_{14}\ne p_{23}$ are satisfied, as in Example~5.1. 

\medskip

\noindent{\bf Example 5.2}  In~\cite{Jon} a just-edge-transitive map $\mathcal M$ of genus $8$, based on the unit cube $\mathcal C$, is constructed with $p_{12}=2$, $p_{23}=p_{24}=p_{34}=3$ and $p_{13}=p_{14}=4$. The vertices of $\mathcal M$ are the eight vertices and twelve edge-midpoints of $\mathcal C$, joined by an edge whenever they are at distance $\sqrt{5}/2$, so  that they have valency $6$ or $4$ respectively. The eight vertices of $\mathcal M$ on each of the six faces of $\mathcal C$ span an octagonal face of $\mathcal M$, and for each of the eight vertices of $\mathcal C$, the six vertices of $\mathcal M$ sharing an edge of $\mathcal C$ with it span a hexagonal face of $\mathcal M$ (see Figure~9). We have ${\rm Aut}\,{\mathcal M}={\rm Aut}\,{\mathcal C}$, and $cr({\mathcal M})=cr({\mathcal C})=2$ corresponding to the fact that the reflections $S_2, S_3$ and $S_4$ are conjugate to each other in $C(p_{ij})$, but are not conjugate to $S_1$.

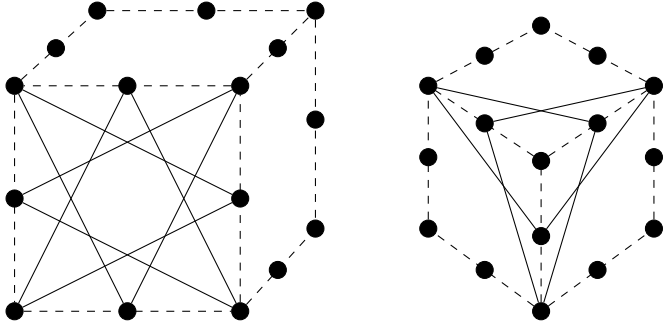
\begin{figure}[h!]
\label{K5}
\begin{center}
 \begin{tikzpicture}[scale=0.5, inner sep=0.8mm]

\node (a) at (-3,0) [shape=circle, draw, fill=black] {}; 
\node (b) at (3,0)  [shape=circle, draw, fill=black] {};
\node (c) at (3,6) [shape=circle, draw, fill=black] {};   
\node (d) at (-3,6)  [shape=circle, draw, fill=black] {};
\node (e) at (-0.8,8)  [shape=circle, draw, fill=black] {};
\node (f) at (5,8)  [shape=circle, draw, fill=black] {};
\node (g) at (5,2.2)  [shape=circle, draw, fill=black] {};

\draw [dashed] (a) to (b);
\draw [dashed] (b) to (c);
\draw [dashed] (c) to (d);
\draw [dashed] (d) to (a);
\draw [dashed] (d) to (e);
\draw [dashed] (c) to (f);
\draw [dashed] (b) to (g);
\draw [dashed] (e) to (f);
\draw [dashed] (f) to (g);

\node (w) at (0,0) [shape=circle, draw, fill=black] {}; 
\node (x) at (3,3)  [shape=circle, draw, fill=black] {};
\node (y) at (0,6) [shape=circle, draw, fill=black] {};   
\node (z) at (-3,3)  [shape=circle, draw, fill=black] {};

\node (h) at (-1.9,7)  [shape=circle, draw, fill=black] {};
\node (i) at (2.1,8)  [shape=circle, draw, fill=black] {};
\node (j) at (4,7)  [shape=circle, draw, fill=black] {};
\node (k) at (4,1.1)  [shape=circle, draw, fill=black] {};
\node (l) at (5,5.1)  [shape=circle, draw, fill=black] {};

\draw (a) to (x);
\draw (x) to (d);
\draw (d) to (w);
\draw (w) to (c);
\draw (c) to (z);
\draw (z) to (b);
\draw (b) to (y);
\draw (y) to (a);


\node (A) at (11,0) [shape=circle, draw, fill=black] {}; 
\node (B) at (11,4)  [shape=circle, draw, fill=black] {};
\node (C) at (8,6) [shape=circle, draw, fill=black] {};   
\node (D) at (14,6)  [shape=circle, draw, fill=black] {};
\node (E) at (8,2.2)  [shape=circle, draw, fill=black] {};
\node (F) at (14,2.2)  [shape=circle, draw, fill=black] {};
\node (G) at (11,7.6)  [shape=circle, draw, fill=black] {};

\draw [dashed] (A) to (B);
\draw [dashed] (B) to (C);
\draw [dashed] (B) to (D);
\draw [dashed] (C) to (G);
\draw [dashed] (D) to (G);
\draw [dashed] (A) to (E);
\draw [dashed] (A) to (F);
\draw [dashed] (E) to (C);
\draw [dashed] (F) to (D);

\node (X) at (11,2) [shape=circle, draw, fill=black] {}; 
\node (Y) at (9.5,5) [shape=circle, draw, fill=black] {};
\node (Z) at (12.5,5) [shape=circle, draw, fill=black] {}; 

\node (H) at (8,4.1)  [shape=circle, draw, fill=black] {};
\node (I) at (9.5,6.8)  [shape=circle, draw, fill=black] {};
\node (I) at (9.5,1.1)  [shape=circle, draw, fill=black] {};
\node (J) at (12.5,6.8)  [shape=circle, draw, fill=black] {};
\node (K) at (12.5,1.1)  [shape=circle, draw, fill=black] {};
\node (L) at (14,4.1)  [shape=circle, draw, fill=black] {};

\draw (A) to (Z);
\draw (Z) to (C);
\draw (C) to (X);
\draw (X) to (D);
\draw (D) to (Y);
\draw (Y) to (A);

\end{tikzpicture}

\end{center}
\caption{An octagonal face and a hexagonal face of $\mathcal M$} 
\end{figure}

\medskip

More generally, by appropriate choices of the parameters $p_{ij}$ one can use similar arguments to show that $cr({\mathcal M})$ can take any of the values $k=1,\ldots, 4$ allowed by Theorem~\ref{edgetrans} for some just-edge-transitive finite map $\mathcal M$.

\subsection{Explicit families of just-edge-transitive maps}

Existence proofs for just-edge-transitive maps based on the residual finiteness of Coxeter groups are non-constructive. Instead, in this section we will give an explicit construction of four infinite families of just-edge-transitive maps $\mathcal M$, with $cr({\mathcal M})=1,\ldots, 4$ respectively.

\begin{thm}\label{justedgetrans}
For each $k=1, 2 ,3$ or $4$ there exist infinitely many compact just-edge-transitive maps ${\mathcal M}_k$ without boundary, such that ${\rm Aut}\,{\mathcal M}_k$ contains exactly $k$ conjugacy classes of reflections. These maps can be chosen to be orientable or non-orientable.
\end{thm}

\noindent{\sl Proof.}  For each $k=1,\ldots, 4$ we will take the maps ${\mathcal M}_k$ to be those corresponding to the kernels $M_k$ of epimorphisms from $N_3$ to finite groups $G_k$, where the images $s_i$ of the generators $S_i$ of $N_3$ lie in $k$ distinct conjugacy classes in $G_k$. 

Let $n$ a prime of the form $4m+1$ for some integer $m\ge 3$ (there are infinitely many such primes, by Dirichlet's Theorem). For each $k=1,\ldots, 4$ we will define four involutions $s_i$ in the symmetric group $S_n$, with $k$ different cycle-structures. Each $s_i$ is a product of $m_i$ disjoint transpositions $t_j:=(j,j+1)$ for $j=1,\ldots, n-1=4m$, which will be assigned to the permutations $s_i$ as follows.

For $k=1, \ldots, 4$ let $\Pi_k$ be the ordered partition $4m=m_1+m_2+m_3+m_4$ of $4m$ defined by:
\[\Pi_1: 4m=m+m+m+m,\]
\[\Pi_2: 4m=(m+1)+(m+1)+(m-1)+(m-1),\]
\[\Pi_3: 4m=(m+2)+m+m+(m-2),\]
\[\Pi_4: 4m=(m+3)+(m+1)+(m-1)+(m-3).\]
Note that in each partition $\Pi_k$ the four summands $m_i$ take $k$ different values, and they all have the same parity.

Given any $k=1,\ldots, 4$, we now define the involutions $s_i$ by assigning the transpositions $t_j=(j, j+1)$ to them as follows:
\begin{itemize}
\item assign $t_1$ to $s_3$;
\item assign $t_2$ to $s_1$ or $s_2$ as $m_1=m_2+2$ or $m_1=m_2$, that is, as $k\in\{1,2 \}$ or $k\in\{3, 4\}$;
\item alternately assign $t_3, \ldots, t_{l+1}$ to $s_3, s_4, s_3,\ldots, s_4, s_3$, where $l=m_3+m_4-1$;
\item alternately assign $t_{l+2}, \ldots, t_{4m-1}=t_{n-2}$ to $s_1, s_2, s_1,\ldots, s_2, s_1$;.
\item assign $t_{4m}=t_{n-1}$ to $s_3$ or $s_4$ as $m_3=m_4+2$ or $m_3=m_4$, that is, as $k\in\{3, 4 \}$ or $k\in\{1, 2\}$.
\end{itemize}

\noindent This is illustrated, for $k\in\{1, 2\}$ and $k\in\{3, 4\}$, in Figures~10 and 11, which show a path ${\mathcal P}_k$ of $n$ vertices labelled $j=1,2, \ldots, n$, with edges $\{j, j+1\}\;(j=1, 2, \ldots, n-1)$ labelled $i=1, 2, 3$ or $4$ as $t_j$ is assigned to $s_i$.

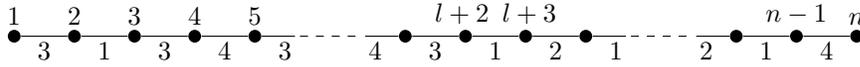
\begin{figure}[h!]
\label{3maps}
\begin{center}
 \begin{tikzpicture}[scale=0.4, inner sep=0.6mm]

\node (1) at (-14,0) [shape=circle, fill=black, label=above: $1$] {};
\node (2) at (-12,0) [shape=circle, fill=black, label=above: $2$] {};
\node (3) at (-10,0) [shape=circle, fill=black, label=above: $3$] {};
\node (4) at (-8,0) [shape=circle, fill=black, label=above: $4$] {};
\node (5) at (-6,0) [shape=circle, fill=black, label=above: $5$] {}; 
\node (6) at (-1,0) [shape=circle, fill=black] {}; 
\node (7) at (1,0) [shape=circle, fill=black, label=above: $l+2\;$] {}; 
\node (8) at (3,0) [shape=circle, fill=black, label=above: $\;l+3$] {}; 
\node (9) at (5,0) [shape=circle, fill=black] {}; 
\node (10) at (10,0) [shape=circle, fill=black] {}; 
\node (11) at (12,0) [shape=circle, fill=black, label=above: $n-1$] {}; 
\node (12) at (14,0) [shape=circle, fill=black, label=above: $n$] {}; 

\node (A) at (-4.5,-0) {};
\node (B) at (-2.5,-0) {};
\node (C) at (6.5,-0) {};
\node (D) at (8.5,-0) {};

\draw (1) to (2);
\draw (2) to (3);
\draw (3) to (4);
\draw (4) to (5);
\draw (5) to (A);
\draw (B) to (6);
\draw (6) to (7);
\draw (7) to (8);
\draw (8) to (9);
\draw (9) to (C);
\draw (D) to (10);
\draw (10) to (11);
\draw (11) to (12);

\draw[dashed] (-4.5,0) to (-2.5,0);
\draw[dashed] (6.5,0) to (8.5,0);

\node at (-13,-0.5) {$3$};
\node at (-11,-0.5) {$1$};
\node at (-9,-0.5) {$3$};
\node at (-7,-0.5) {$4$};
\node at (-5,-0.5) {$3$};
\node at (-2,-0.5) {$4$};
\node at (0,-0.5) {$3$};
\node at (2,-0.5) {$1$};
\node at (4,-0.5) {$2$};
\node at (6,-0.5) {$1$};
\node at (9,-0.5) {$2$};
\node at (11,-0.5) {$1$};
\node at (13,-0.5) {$4$};

\end{tikzpicture}

\end{center}
\caption{The labelled path ${\mathcal P}_k$ where $k=1$ or $2$} 
\end{figure}

\begin{figure}[h!]
\label{3maps}
\begin{center}
 \begin{tikzpicture}[scale=0.4, inner sep=0.6mm]

\node (1) at (-14,0) [shape=circle, fill=black, label=above: $1$] {};
\node (2) at (-12,0) [shape=circle, fill=black, label=above: $2$] {};
\node (3) at (-10,0) [shape=circle, fill=black, label=above: $3$] {};
\node (4) at (-8,0) [shape=circle, fill=black, label=above: $4$] {};
\node (5) at (-6,0) [shape=circle, fill=black, label=above: $5$] {}; 
\node (6) at (-1,0) [shape=circle, fill=black] {}; 
\node (7) at (1,0) [shape=circle, fill=black, label=above: $l+2\;$] {}; 
\node (8) at (3,0) [shape=circle, fill=black, label=above: $\;l+3$] {}; 
\node (9) at (5,0) [shape=circle, fill=black] {}; 
\node (10) at (10,0) [shape=circle, fill=black] {}; 
\node (11) at (12,0) [shape=circle, fill=black, label=above: $n-1$] {}; 
\node (12) at (14,0) [shape=circle, fill=black, label=above: $n$] {}; 

\node (A) at (-4.5,-0) {};
\node (B) at (-2.5,-0) {};
\node (C) at (6.5,-0) {};
\node (D) at (8.5,-0) {};

\draw (1) to (2);
\draw (2) to (3);
\draw (3) to (4);
\draw (4) to (5);
\draw (5) to (A);
\draw (B) to (6);
\draw (6) to (7);
\draw (7) to (8);
\draw (8) to (9);
\draw (9) to (C);
\draw (D) to (10);
\draw (10) to (11);
\draw (11) to (12);

\draw[dashed] (-4.5,0) to (-2.5,0);
\draw[dashed] (6.5,0) to (8.5,0);

\node at (-13,-0.5) {$3$};
\node at (-11,-0.5) {$2$};
\node at (-9,-0.5) {$3$};
\node at (-7,-0.5) {$4$};
\node at (-5,-0.5) {$3$};
\node at (-2,-0.5) {$4$};
\node at (0,-0.5) {$3$};
\node at (2,-0.5) {$1$};
\node at (4,-0.5) {$2$};
\node at (6,-0.5) {$1$};
\node at (9,-0.5) {$2$};
\node at (11,-0.5) {$1$};
\node at (13,-0.5) {$3$};

\end{tikzpicture}

\end{center}
\caption{The labelled path ${\mathcal P}_k$ where $k=3$ or $4$} 
\end{figure}
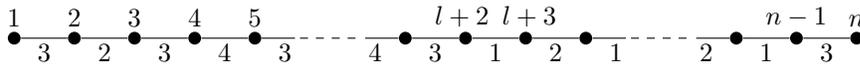

For any $k$, each $s_i$ is a product of $m_i$ disjoint transpositions $t_j$, so it is an involution. The graph ${\mathcal P}_k$ is connected, so the permutation group $G_k:=\langle s_1,\ldots, s_4\rangle$ is transitive, and is therefore primitive since $n$ is prime. Since 
\[s_3s_4=(1, 2)(3, 4, \ldots, l+2)(n-1, n)\]
with $l=m_3+m_4-1$ odd, $(s_3s_4)^2$ is a cycle of length $l$ with $n-l\ge 3$ fixed points. A theorem of Jordan (see~\cite[Theorem~3.3E]{DM} or~\cite[Theorem~13.9]{Wie}) states that if a finite primitive permutation group contains a cycle of prime length $l$ with at least three fixed points, then it contains the alternating group. This result has been extended in~\cite{Jon2014} to cycles of all lengths $l\ge 2$, so $G_k$ contains $A_n$. Since $G_k=\langle s_1,\ldots, s_4\rangle$, we have $G_k=A_n$ or $S_n$ as the summands $m_i$, in $\Pi_k$, which have the same parity,  are all even or all odd. Since the permutations $s_i$ have $k$ different cycle-structures, they lie in $k$ different conjugacy classes in $G_k$. (Here we use the fact that involutions in $A_n$ are conjugate in $A_n$ if and only if they are conjugate in $S_n$, that is, if and only if they have the same cycle structure.)

Let $M_k$ be the kernel of the epimorphism $N_3\to G_k$ given by $S_i\mapsto s_i$ for $i=1,\ldots, 4$, so that $N_3\le N_{\Gamma}(M_k)$. In order that $N_3=N_{\Gamma}(M_k)$ we require that no double transposition of the generators $s_i$ extends to an automorphism of $G_k$. Since $n>6$ we have ${\rm Aut}\,G_k=S_n$, acting by conjugation on $G_k=A_n$ or $S_n$, so we want no permutation $g\in S_n$, acting by conjugation, to induce a double transposition on the generators $s_i$. Any such permutation $g$ would have to extend to an automorphism of the graph ${\mathcal P}_k$, inducing a double transposition on the edge-labels. It is clear from Figures~10 and 11 that no such automorphism exists, so $N_3=N_{\Gamma}(M_k)$ as required. The map ${\mathcal M}_k$ with map subgroup $M_k$ is therefore just-edge-transitive. It is a regular covering of the map ${\mathcal N}_3$  in Figure~5, with covering group and automorphism group ${\rm Aut}\,{\mathcal M}_k\cong N_3/M_k\cong G_k\cong A_n$ or $S_n$. It has $k$ conjugacy classes of reflections, corresponding to the classes in $G_k$ containing the involutions $s_i$. It ls compact and without boundary, since $M_k$ has finite index in $\Gamma$ and contains no conjugate of any reflection $R_i$. It is orientable or non-orientable as the permutations $s_i$ (or equivalently the summands $m_i$) are all odd or all even; for each $k$, by appropriately choosing primes $n\equiv 1$ or $n\equiv 5\pmod{8}$ one can find infinitely many examples in either case. \hfill$\square$

\medskip

\noindent{\bf Remarks}

\smallskip

\noindent{\bf 1} The extension of Jordan's Theorem mentioned above depends, through the classification of doubly transitive permutation groups, on the classification of finite simple groups (see~\cite{Jon2014}). This deep result could be avoided by using a similar but more complicated construction which ensures that the primitive group $G_k$ contains a cycle of {\sl prime\/} length with at least three fixed points, so that $G_k\ge A_n$ by Jordan's original theorem.

\smallskip

\noindent{\bf 2} In this construction, the orders $p_{ij}$ of the various products $s_is_j$, and hence the valencies of the vertices, faces and Petrie polygons, can easily be read off from the corresponding graph ${\mathcal P}_k$. For instance, if $k\in\{1, 2\}$ then $s_1s_4$ and $s_2s_3$ have orders $6$ and $2$, so the faces are $12$-gons and $4$-gons, whereas if $k\in\{3, 4\}$ they are $4$-gons and $8$-gons.  

\smallskip

\noindent{\bf 3} Together with the orientably regular chiral maps, those with $N=\Gamma^+$, the examples constructed here show that if $\mathcal M$ is an edge-transitive map then $cr({\mathcal M})$ can take any value $0,\ldots, 4$; by Theorem~\ref{edgetrans}, no other values are possible. 

\smallskip

\noindent{\bf 4} In general one should not expect to construct edge-transitive maps $\mathcal M$ with conjugacy classes of reflections of arbitrary given sizes, as was achieved for vertex-transitive maps in Theorem~\ref{vtrans}. For instance, if we try to obtain $k=4$ classes with sizes $c_1=c_2=c_3=1$ and $c_4>1$, then ${\rm Aut}\,{\mathcal M}$ is generated by four involutions with only three of them in the centre, which is clearly impossible.

\smallskip

\noindent{\bf 5} This section has concentrated on the just-edge-transitive maps, since they include all the edge-transitive maps attaining the upper bound $cr({\mathcal M})\le 4$ of Theorem~\ref{edgetrans}. In principle it would be straightforward to apply a similar analysis to the other thirteen cases of the Graver-Watkins classification listed in Section~5.2: the first two cases, types~1 and $2^P$ex, consisting of the regular maps and the orientably regular chiral maps, have already been intensively studied by many authors (see~\cite{Sir} for an excellent recent survey), while Sir\'a\v n, Tucker and Watkins have given explicit constructions of infinite families of compact maps of all fourteen types in~\cite{STW}.

\section{Reflections of hypermaps}

Hypermaps are generalisations of maps, in which edges (now called hyperedges) are allowed to be incident with any number of vertices and faces, rather than at most two. They give useful combinatorial representations of triangle groups, and in the compact, orientable case they play a major role in Grothendieck's theory of dessins d'enfants~\cite{JS}. Much of the discussion here of reflections of maps carries over to hypermaps with only minimal changes: for instance, the group $\Gamma\cong V_4*C_2$ is replaced with the group $\Delta\cong C_2*C_2*C_2$ obtained from $\Gamma$ by omitting the relation $(R_0R_2)^2=1$. All maps are hypermaps, so the constructions in Theorem~\ref{vtrans}, Corollary~\ref{vtranscor} and Theorem~\ref{justedgetrans} can be regarded as proving the existence of hypermaps with given classes of reflections.

The upper bound $cr({\mathcal M})\le 3$ for regular maps is also valid for regular hypermaps. However, since $E=\langle R_0, R_2\rangle$ is now an {\sl infinite\/} dihedral group, the edge-transitivity condition $\Delta=NE$ is much weaker than for maps, and Theorem~\ref{edgetrans} does not extend to hypermaps. Indeed, the existence of a group of automorphisms of $\Delta$, acting as $S_3$ on the generators $R_i$, shows that vertex-, edge- and face-transitivity are all essentially equivalent to each other, modulo certain obvious duality and triality operations on hypermaps, so it is the analogues of Theorem~\ref{vtrans} and Corollary~\ref{vtranscor} which apply.

\section{Riemann surfaces and algebraic curves}

Any orientable map $\mathcal M$ without boundary induces a complex structure on its underlying surface $S$, making it a Riemann surface. This can be done by regarding $\Gamma$ as the automorphism group of the universal map ${\mathcal M}_{\infty}$ (see Figure~4), and taking $\mathcal M$ to be the quotient ${\mathcal M}_{\infty}/M$ of ${\mathcal M}_{\infty}$ by the map subgroup $M\le\Gamma$ fixing some $\phi\in\Phi$. This construction gives $S$ a complex structure, inherited from $\mathbb H$, so that it is a Riemann surface, and each automorphism of $\mathcal M$ induces a conformal or anticonformal automorphism of $S$ as it preserves or reverses the orientation. (If $\mathcal M$ is non-orientable, or has non-empty boundary, then a similar process induces the structure of a Klein surface on $S$.) 

A compact Riemann surface $S$ can be regarded as a complex algebraic curve. Conjugacy classes of anticonformal involutions (often called symmetries) of $S$ correspond to isomorphism classes of real forms of this curve; those with fixed points (called reflections) correspond to real forms with real points (see~\cite{BCGG}, for example). In~\cite{Nat}, Natanzon showed that a compact Riemann surface $S$ of genus $g\ge 2$ has at most $2\sqrt g+2$ conjugacy classes of reflections, and that this bound is attained for infinitely many values of $g$, of the form $(2^n-1)^2$; see~\cite{BGS} for an alternative proof and for other related results. Subsequently Bujalance, Gromadzki and Izquierdo~\cite{BGI} showed that if $g=2^{r-1}u+1\ge 2$ with $u$ odd then $S$ has at most $2^{r+1}$ conjugacy classes of reflections, attained if and only if $u\ge 2^{r+1}-3$ (see also~\cite[Theorem 2.2.1]{BCGG}); this extended an earlier upper bound of four conjugacy classes for even $g$, given by Gromadzki and Izquerdo~\cite{GI}.

These bounds apply to {\sl all\/} compact Riemann surfaces of a given genus. In the cases considered in this paper, the complex structures are obtained from maps, and by Bely\u\i's Theorem~\cite{Bel} this is equivalent to the corresponding curves being defined over the field $\overline{\mathbb Q}$ of algebraic numbers; conjugacy classes of reflections then correspond to real forms defined over  the field ${\mathbb R}\cap\overline{\mathbb Q}$ of real algebraic numbers, by the extensions of Bely\u\i's Theorem due to K\"ock and Singerman~\cite{KS} and to K\"ock and Lau~\cite{KL}. The results obtained here apply in this context, rather than the more general context of compact Riemann or Klein surfaces.

In a fundamental paper, Bujalance and Singerman~\cite{BS} have given a detailed study of the reflections of a Riemann surface, considering the possible {\sl symmetry types}, or lists of {\sl species }$\pm k$ which can arise: here $k$ is the number of fixed curves of a reflection representing each conjugacy class, and the symbol $\pm$ denotes whether or not the quotient surface is orientable. It would be interesting to combine their approach with that used in the present paper, by considering the possible symmetry types of maps under various transitivity conditions on flags, vertices, edges, etc. Similarly, it would be interesting to extend to this wider context the results of Meleko\u glu and Singerman~\cite{Mel, MS, MS14} on patterns of reflections of regular maps, describing how their fixed curves pass through sequences of vertices, edges and faces. Finally, one might consider whether the methods developed here could be applied to other categories, where objects such as abstract polytopes can also admit reflections.

\bigskip

\noindent{\bf Acknowledgement} The author is grateful to the organisers of the conference SIGMAP 2014, to Adnan Meleko\u glu, whose talk at that meeting did much to motivate this work, to the organisers of the conference on Graph Embeddings (St Petersburg, November 2014) for the opportunity to give a talk on which this paper is based, and to Ian Leary and David Singerman for some very helpful comments.


\newpage

\begin{thebibliography}{99}

\bibitem{Bel} G.~V.~Bely\u\i, On Galois extensions of a maximal cyclotomic field, {\sl Izv.~Akad. Nauk SSSR Ser Mat.} 43 (1979), 267--276, 479. 


\bibitem{BS} R.~P.~Bryant and D.~Singerman, Foundations of the theory of maps on surfaces with boundary, {\sl Q.~J.~Math.} (2) 36 (1985), 17--41.

\bibitem{BCGG} E.~Bujalance, F.~J.~Cirre, J.~M.~Gamboa and G.~Gromadzki, {\sl Symmetries of Compact Riemann Surfaces}, Lecture Notes in Mathematics 2007, Springer, 2010. 

\bibitem{BGI} E.~Bujalance, G.~Gromadzki and M.~Izquierdo, On real forms of a complex algebraic curve, {\sl J.~Aust.~Math.~Soc.} 70 (2001), 134--142.

\bibitem{BGS} E.~Bujalance, G.~Gromadzki and D.~Singerman, On the number of real curves associated to a complex algebraic curve, {\sl Proc.~Amer.~Math.~Soc.} 120 (1994), 507--513.

\bibitem{BS} E.~Bujalance and D.~Singerman, The symmetry type of a Riemann surface, {\sl Proc. Lond.~Math.~Soc.} (3) 51 (1985), 501--519.

\bibitem{Bur} W.~Burnside, On groups of order $p^aq^b$, {\sl Proc.~London Math.~Soc.} 2 (1904), 388--392.





\bibitem{CM}  H.~S.~M.~Coxeter and W.~O.~J.~Moser, {\sl Generators and Relations for Discrete Groups}, 4th ed., Springer-Verlag, Berlin -- Heidelberg -- New York, 1980.  %


\bibitem{DM} J.~D.~Dixon and B.~Mortimer, {\sl Permutation Groups}, Graduate Texts in Mathematics 163, Springer, 1996.




\bibitem{GW} J.~E.~Graver and M.~E.~Watkins, Locally finite, planar, edge-transitive graphs, {\sl Mem.~Amer.~Math.~Soc.} 126 (1997), no.~601.

\bibitem{GI} G.~Gromadzki and M.~Izquierdo, Real forms of a Riemann surface of even genus, {\sl Proc.~Amer.~Math.~Soc.} 126 (1998), 3475--3479.

\bibitem{Gro} A.~Grothendieck, Esquisse d'un Programme, in {\sl Geometric Galois Actions~I, Around Grothendieck's Esquisse d'un Programme} (ed.~P.~Lochak and L.~Schneps), London Math.~Soc.~Lecture Note Ser.~242 (Cambridge University Press, Cambridge, 1997), pp.~5--48. 

\bibitem{Hup} B.~Huppert, {\sl Endliche Gruppen I}, Springer-Verlag, Berlin, 1967.


\bibitem{Jon} G.~A.~Jones, Just-edge-transitive maps and Coxeter groups, {\sl Ars Combin.} 16-B (1983), 139--150.

\bibitem{Jon2014} G.~A.~Jones, Primitive permutation groups containing a cycle, {\sl Bull.~Aust. Math.~Soc.} 89 (2014), 159--165.

\bibitem{JS} G.~A.~Jones and D.~Singerman, Maps, hypermaps and triangle groups, in {\sl The Grothendieck Theory of Dessins d'Enfants} (ed.~L.~Schneps), London Math.~Soc.~Lecture Note Ser.~200 (Cambridge University Press, Cambridge, 1994), pp.~115--145. 


\bibitem{KL} B.~K\"ock and E.~Lau, A note on Belyi's theorem for Klein surfaces, {\sl Q.~J.~Math.} 61 (2010), 103--107.

\bibitem{KS} B.~K\"ock and D.~Singerman, Real Belyi theory, {\sl Q.~J.~Math.} 58 (2007), 463--478.


\bibitem{LS} R.~C.~Lyndon and P.~E.~Schupp, {\sl Combinatorial Group Theory}, Springer-Verlag, Berlin -- Heidelberg -- New York, 1977. 



\bibitem{MKS} W.~Magnus, A.~Karrass and D.~Solitar, {\sl Combinatorial Group Theory}, Dover, New York, 1976. 

\bibitem{Mal} A.~I.~Mal'cev, On isomorphic matrix representations of infinite groups of matrices (Russian), {\sl Mat.~Sb.} 8 (1940), 405--422; {\sl Amer.~Math.~Soc.~Transl.} (2) 45 (1965), 1--18.

\bibitem{Mel} A.~Meleko\u glu, A geometric approach to the reflections of regular maps, {\sl Ars Combin.} 89 (2008), 355--367.

\bibitem{MS} A.~Meleko\u glu and D.~Singerman, Reflections of regular maps and Riemann surfaces, {\sl Rev.~Mat.~Iberoam.} 24 (2006), 921--939.

\bibitem{MS14} A.~Meleko\u glu and D.~Singerman, The structure of mirrors on Platonic surfaces, submitted.

\bibitem{Nat} S.~M.~Natanzon, The order of a finite group of homeomorphisms of a surface onto itself, and real forms of a complex algebraic curve, (Russian), {\sl Dokl.~Akad.~Nauk SSSR} 242 (1978), 765--768. English translation: {\sl Soviet Math.~Dokl.} 19 (1978), 1195--1199 (1979)

\bibitem{Rem} V.~N.~Remeslennikov, Finite approximability of groups with respect to conjugacy (Russian), {\sl Sibirsk.~Mat.~\v Z.} 12 (1971), 1085--1099.






\bibitem{Sir} J.~\v Sir\' a\v n, How symmetric can maps on surfaces be?, in {\sl Surveys in Combinatorics 2013} (eds.~S.~R.~Blackburn, S.~Gerke and M.~Wildon), London Math.~Soc.~Lecture Note Ser.~409 (Cambridge University Press, Cambridge, 2013), pp.~161--238.

\bibitem{STW} J.~\v Sir\'a\v n, T.~W.~Tucker and M.~E.~Watkins, Realizing finite edge-transitive orientable maps, {\sl J.~Graph Theory} 37 (2001), 1--34.

\bibitem{Ste} P.~Stebe, A residual property of certain groups, {\sl Proc.~Amer.~Math.~Soc.} 26 (1970), 37--42.



\bibitem{Wie} H.~Wielandt, {\sl Finite Permutation Groups}, Academic Press, New York, 1964.


\end{thebibliography}
\end{document}